\documentclass[%
]{mpi2015-cscpreprint}

\usepackage{amsfonts, amsmath, amsopn, amssymb, bm, mathrsfs}
\usepackage{algorithm, algorithmicx, algpseudocode}
\usepackage{amsthm}
	\newtheorem{definition}{Definition}
	\newtheorem{theorem}{Theorem}
	\newtheorem{remark}{Remark}
\usepackage{datetime2}
\usepackage[inline, shortlabels]{enumitem}
\usepackage{graphicx}
\usepackage{booktabs, makecell} 
\newcolumntype{P}[1]{>{\centering\arraybackslash}p{#1}}
\newcolumntype{M}[1]{>{\centering\arraybackslash}m{#1}}
\usepackage[section]{placeins}
\usepackage{xcolor}
\usepackage{xspace}

\DeclareMathOperator{\spn}{span}	

\DeclareMathOperator{\spnS}{\spn^{\mathbb{S}}}

\DeclareMathOperator{\trace}{trace}

\newcommand{\chol}{\texttt{chol}\xspace}

\newcommand{\nan}{\texttt{NaN}\xspace}

\newcommand{\spC}{\mathbb{C}}				

\newcommand{\spCnn}{\spC^{n \times n}}

\newcommand{\spCns}{\spC^{n \times s}}

\newcommand{\spCss}{\spC^{s \times s}}

\newcommand{\spK}{\mathscr{K}}				
\newcommand{\spS}{\mathbb{S}}

\newcommand{\spKS}{\spK^{\spS}}        

\newcommand{\ve}{\bm{e}}
\newcommand{\vehat}{\widehat{\ve}}

\newcommand{\vB}{\bm{B}}

\newcommand{\vD}{\bm{D}}

\newcommand{\vE}{\bm{E}}
\newcommand{\vEhat}{\widehat{\bm{E}}}

\newcommand{\vJ}{\bm{J}}

\newcommand{\vM}{\bm{M}}

\newcommand{\vQ}{\bm{Q}}

\newcommand{\vR}{\bm{R}}

\newcommand{\vU}{\bm{U}}

\newcommand{\vV}{\bm{V}}

\newcommand{\vW}{\bm{W}}

\newcommand{\vX}{\bm{X}}

\newcommand{\vY}{\bm{Y}}

\newcommand{\vZ}{\bm{Z}}

\newcommand{\vXi}{\bm{\Xi}}

\newcommand{\Beta}{B}

\newcommand{\HH}{\mathcal{H}}

\newcommand{\MM}{\mathcal{M}}

\newcommand{\RR}{\mathcal{R}}

\newcommand{\TT}{\mathcal{T}}

\newcommand{\bQQ}{\bm{\mathcal{Q}}}

\newcommand{\bVV}{\bm{\mathcal{V}}}

\newcommand{\bXX}{\bm{\mathcal{X}}}

\newcommand{\inv}{{-1}}

\newcommand{\cinv}{{-*}}

\newcommand{\bigO}[1]{\mathcal{O}\left(#1\right)}

\newcommand{\norm}[1]{\left\lVert#1\right\rVert}

\newcommand{\normF}[1]{\norm{#1}_{\text{F}}}

\makeatletter
\newsavebox{\@brx}
\newcommand{\llangle}[1][]{\savebox{\@brx}{\(\m@th{#1\langle}\)}%
	\mathopen{\copy\@brx\kern-0.5\wd\@brx\usebox{\@brx}}}
\newcommand{\rrangle}[1][]{\savebox{\@brx}{\(\m@th{#1\rangle}\)}%
	\mathclose{\copy\@brx\kern-0.5\wd\@brx\usebox{\@brx}}}
\makeatother

\newcommand{\IPS}[2]{{\llangle #1,#2 \rrangle}_{\spS}}
\newcommand{\IPSnoarg}{\IPS{\cdot}{\cdot}}
\newcommand{\N}{N}

\newcommand{\BCGS}{\texttt{BCGS}\xspace}	
\newcommand{\BCGSPIP}{\texttt{BCGS-PIP}\xspace}	
\newcommand{\BCGSPIO}{\texttt{BCGS-PIO}\xspace}	
\newcommand{\BCGSIROLS}{\texttt{BCGSI+LS}\xspace}	
\newcommand{\BMGS}{\texttt{BMGS}\xspace}	
\newcommand{\BMGSSVL}{\texttt{BMGS-SVL}\xspace}	
\newcommand{\BMGSLTS}{\texttt{BMGS-LTS}\xspace}	
\newcommand{\BMGSCWY}{\texttt{BMGS-CWY}\xspace}	
\newcommand{\BMGSICWY}{\texttt{BMGS-ICWY}\xspace}	

\newcommand{\IO}[1]{\texttt{IntraOrtho}\left(#1\right)}	
\newcommand{\IOnoarg}{\texttt{IntraOrtho}\xspace}	
\newcommand{\CGSIROLS}{\texttt{CGSI+LS}\xspace}	
\newcommand{\MGSSVL}{\texttt{MGS-SVL}\xspace}	
\newcommand{\MGSLTS}{\texttt{MGS-LTS}\xspace}	
\newcommand{\HouseQR}{\texttt{HouseQR}\xspace}	
\newcommand{\AllReduceQR}{\texttt{AllReduceQR}\xspace}
\newcommand{\TSQR}{\texttt{TSQR}\xspace}	
\newcommand{\CholQR}{\texttt{CholQR}\xspace}	

\newcommand{\cl}{\texttt{cl}}
\newcommand{\gl}{\texttt{gl}}

\newcommand{\eps}{\varepsilon}

\newcommand{\MATLAB}{\textsc{Matlab}\xspace}

\newcommand{\LSBA}{\texttt{LowSyncBlockArnoldi}\xspace}
\newcommand{\tridiag}{\texttt{tridiag}\xspace}
\newcommand{\powersystem}{\texttt{1138\_bus}\xspace}
\newcommand{\lapl}{\texttt{lapl\_2d}\xspace}

\newcommand{\matvec}{\texttt{matvec}\xspace}
\newcommand{\innerprod}{\texttt{inner\_prod}\xspace}
\newcommand{\intraortho}{\texttt{intra\_ortho}\xspace}
\newcommand{\basiseval}{\texttt{basis\_eval}\xspace}

\newcommand{\fom}{\texttt{FOM}\xspace}
\newcommand{\gmres}{\texttt{GMRES}\xspace}

\definecolor{ForestGreen}{RGB}{34,139,34}
\newcommand{\updated}[1]{\textcolor{black}{#1}}

\begin{document}
  

\title{Adaptively restarted block Krylov subspace methods with low-synchronization skeletons}
  
\author[$\ast$]{Kathryn Lund}
\affil[$\ast$]{Computational Methods in Systems and Control Theory, Max Planck Institute for Dynamics of Complex Technical Systems, Sandtorstr.~1, Magdeburg 39106, Germany\authorcr
  \email{lund@mpi-magdeburg.mpg.de},
  \orcid{0000-0001-9851-6061}
}
  
\shorttitle{Low-sync block Arnoldi}
\shortauthor{K.~Lund}
\shortdate{}
  
\keywords{Gram-Schmidt, Krylov subspace methods, Arnoldi method, block methods, stability, loss of orthogonality, low-synchronization methods, high-performance computing, communication-avoiding methods}

\msc{65F10, 65F25, 65F50, 15-04}
  
\abstract{With the recent realization of exascale performace by Oak Ridge National Laboratory's Frontier supercomputer, reducing communication in kernels like QR factorization has become even more imperative.  Low-synchronization Gram-Schmidt methods, first introduced in [\textsc{K.~{\'{S}}wirydowicz, J.~Langou, S.~Ananthan, U.~Yang, and S.~Thomas}, {\em {Low Synchronization Gram-Schmidt and Generalized Minimum Residual Algorithms}}, Numer.~Lin.~Alg.~Appl., Vol.~28(2), e2343, 2020)], have been shown to improve the scalability of the Arnoldi method in high-performance distributed computing.  Block versions of low-synchronization Gram-Schmidt show further potential for speeding up algorithms, as column-batching allows for maximizing cache usage with matrix-matrix operations.  In this work, low-synchronization block Gram-Schmidt variants from [\textsc{E.~Carson, K.~Lund, M.~Rozlo\v{z}n\'{i}k, and S.~Thomas}, {\em Block {Gram}-{Schmidt} algorithms and their stability properties}, Lin.~Alg.~Appl., 638, pp.~150--195, 2022] are transformed into block Arnoldi variants for use in block full orthogonalization methods (BFOM) and block generalized minimal residual methods (BGMRES).  An adaptive restarting heuristic is developed to handle instabilities that arise with the increasing condition number of the Krylov basis.  The performance, accuracy, and stability of these methods are assessed via a flexible benchmarking tool written in MATLAB.  The modularity of the tool additionally permits generalized block inner products, like the global inner product.
}

\novelty{We introduce new block Arnoldi algorithms for communication-avoiding Krylov subspace methods, along with a software tool and a proposed set of benchmarking measures.}

\maketitle

  
\section{Introduction and motivation} \label{sec:intro}
Oak Ridge National Laboratory reported in May 2022 that its Frontier supercomputer is the first machine to have achieved true exascale performance\footnote{\url{http://www.top500.org/news/ornls-frontier-first-to-break-the-exaflop-ceiling/}. Accessed 8 August 2022.}.  That is, for the first time ever, a supercomputer performed more than 1 exaflop (i.e., $10^{18}$ double-precision floating-point operations) in a single second.  This astounding development is clear motivation for our work.  Exascale computing is no longer a next-generation dream; it is reality, and the need for highly parallelized algorithms that take full advantage of exaflop computational potential while reducing global communication between nodes is urgent.

To this end we build on the low-synchronization (``low-sync") Gram-Schmidt methods of Barlow \cite{Bar19}, \`{S}wirydowicz et al.~\cite{SwiLAetal20}, Yamazaki et al.~\cite{YamTHetal20}, Thomas et al.~\cite{ThoCRetal22}, and Bielich et al.~\cite{BieLTetal22}, as well as our own earlier work with block versions of these methods \cite{CarLR21, CarLRetal22}.  Gram-Schmidt methods are an essential backbone in orthogonalization routines like QR factorization and in iterative methods like Krylov subspace methods for linear systems, matrix functions, and matrix equations \cite{Saa03a, Gue13, Sim16}.  Block Krylov subspace methods in particular make better use of L3 cache via matrix-matrix operations and feature often in communication-avoiding Krylov subspaces, such as $s$-step \cite{Car15, Hoe10}, enlarged methods \cite{GriMN16}, and randomized methods \cite{BalGr21}.

As in most realms of life, there is no such thing as a free lunch here.  While low-sync variations have the potential to speed up highly parallelized implementations of Gram-Schmidt \cite{YamTHetal20}, they introduce new floating-point errors and thus potential instability, due to the reformulation of inner products and normalizations.  Instability surfaces in the loss of orthogonality between basis vectors and can lead to breakdowns or wildly inaccurate approximations in downstream applications \cite{Hig02, HucNe19}.  Stability bounds for some low-sync variants have been established, but it often takes much longer to carry out a rigorous stability analysis than to derive and deploy new methods \cite{Bar19, CarLR21, CarLRetal22, ThoCRetal22}.  It can also happen that a backward error bound is established  and later challenged by an obscure edge case \cite{GirLRetal05, SmoBL06}.  With this tension in mind, we have not only extended low-sync variants of block Gram-Schmidt to block Arnoldi but also developed a benchmarking tool for the community to explore the efficiency, stability, and accuracy of these new algorithms, in a similar vein as the \texttt{BlockStab}\footnote{\url{https://github.com/katlund/BlockStab}} comparison tool developed in tandem with a recent block Gram-Schmidt survey \cite{CarLRetal22}.  We refer to this new tool as \LSBA\footnote{\url{https://gitlab.mpi-magdeburg.mpg.de/lund/low-sync-block-arnoldi}} and encourage the reader to explore the tool in parallel with the text.

Established in this earlier work is the fact that block variants of low-sync Gram-Schmidt are less stable than their column-wise counterparts.  However, when these skeletons are transferred to block Arnoldi and used to solve linear systems, we gain the option to restart the process.  Restarting can \updated{be effective} at mitigating stability issues in communication-avoiding algorithms \cite{Car18, Car20}.  As long as each node redundantly computes residual or error estimates and checks the stability via local quantities, restarting does not introduce additional synchronization points.  Furthermore, adaptive restarting allows for robustness, as we can use basic look-ahead heuristics to foresee a breakdown and salvage progress without giving up completely at the first sign of trouble.

Given the modularity of our framework, we are also able to treat generalized block inner products, as described in \cite{FroLS17, FroLS20}.  We focus in particular on the classical and global inner products.

The paper is organized as follows.  In Section~\ref{sec:background} we summarize terms, definitions, and concepts from high-performance (HPC) computing, generalized block inner products, block Gram-Schmidt algorithms, and block Krylov subspace methods with static restarting.  We present new low-synchronization block Arnoldi skeletons in Section~\ref{sec:low_sync}, and derive an adaptive restarting heuristic in Section~\ref{sec:stab_repro}.  Section~\ref{sec:examples} features a more in-depth discussion of the \LSBA benchmarking tool as well as examples demonstrating how to compare different block Arnoldi variants.  We summarize our findings in Section~\ref{sec:conclusion}.


\section{Background} \label{sec:background}
This work is \updated{a combination of} the generalized inner product framework of Frommer, Lund, and Szyld \cite{FroLS17, FroLS20} and the skeleton-muscle framework for block Gram-Schmidt (BGS) by Carson, Lund, Rozlo\v{z}n\'{i}k, and Thomas \cite{CarLR21, CarLRetal22}.  Throughout the text, we focus on solving a linear system with multiple right-hand sides
\begin{equation} \label{eq:AX=B}
	A \vX = \vB,
\end{equation}
where $A \in \spCnn$ is large and sparse (i.e., with $\bigO{n}$ nonzero entries) and $\vB \in \spCns$ is \updated{a} tall-skinny (i.e., $s \ll n$) matrix.

We employ standard numerical linear algebra notation throughout.  In particular, $A^*$ denotes the Hermitian transpose of $A$, $\norm{\cdot}$ refers to the Euclidean 2-norm, unless otherwise specified, and $\vehat_k$ denotes the $k$th standard unit vector with the $k$th entry equal to $1$ and all others $0$.

In the following subsections, we define key concepts in HPC, block Gram-Schmidt methods, and block Krylov subspace methods.

\subsection{Communication in high-performance computing} \label{sec:comm_hpc}
As floating-point operations have become faster and less energy-intensive, {\em communication}-- the memory operations \updated{between} levels of cache on a node or between parallelized processors on a network-- has become a bottleneck in distributed computing.  How expensive a memory operation is depends on the physical aspects of a specific system, specifically the {\em latency}, or the amount of time needed to pack and transmit a message, and the {\em bandwidth}, or how much information can be transmitted at a time.  To improve algorithm performance in bandwidth-limited algorithms like Krylov subspace methods, it is therefore advantageous to increase the {\em computational intensity}, or the ratio between floating-point and memory operations \cite{BalCDetal14}.  We pay particular attention to {\em synchronization points} (``sync points"), i.e., the steps in an algorithm that initiate a broadcast or reduce pattern to synchronize a quantity on all processors.   Reducing calls to kernels with sync points is a straightforward way to improve computational intensity \cite{AnzBFetal20}.

Sync points in Krylov subspace methods arise primarily in the orthonormalization procedure, such as Arnoldi or Lanczos, both of which are reformulations of the Gram-Schmidt method, a standard method for orthonormalizing a basis one (block) vector at a time.  For large $n$, vectors are typically partitioned row-wise and distributed among processors, meaning that any time an operation like an inner product or normalization is performed-- which is at least once per (block) vector in Gram-Schmidt-- a sync point is inevitable.

Other possibly communication-intensive kernels include applications of the operator $A$\footnote{The term \texttt{matvec} is often used to refer to the multiplication of $A$ with a vector.  Because we will be focusing on block vectors, we refrain from this term to avoid confusion.} and applications of $\bVV_m$, an $n \times ms$ Krylov basis matrix.  We count each operation separately from sync points (block inner products and vector norms) in \LSBA; see Section~\ref{sec:examples}.

\subsection{Generalized block inner products} \label{sec:block_inner_prod}
A {\em block vector} is a tall-skinny matrix $\vX \in \spC^{n \times s}$, and a {\em block matrix} is a matrix of $s \times s$ matrices, e.g.,
\[
\HH =
\begin{bmatrix}
	H_{1,1}	& H_{1,2}	&	\cdots	& H_{1,p} \\
	H_{2,1}	& H_{2,2}	&	\cdots	& H_{2,p} \\
	\vdots	&	\vdots	&	\ddots	& \vdots \\
	H_{q,1}	& H_{q,2}	&	\cdots	& H_{q,p}
\end{bmatrix}
\in \spC^{qs \times ps}.
\]
We use a mixture of \MATLAB- and block-indexing notation to handle block objects.  In particular, we write $\bVV_k$ to denote the first $k$ block vectors of the block-partitioned matrix $\bVV = \begin{bmatrix} \vV_1 & \vV_2 & \cdots & \vV_m \end{bmatrix}$ instead of \updated{$\bVV_{:,1:ks}$} (i.e., the first $ks$ columns).  In a similar vein, $s \times s$ block entries of $\HH$ are denoted as $H_{j,k}$ instead of as $H_{(j-1)s+1:js,(k-1)s+1:ks}$.  We denote block generalizations of the standard unit vectors $\vehat_k$ as $\vEhat_k := \vehat_k \otimes I_s$, where $\otimes$ is the Kronecker product and $I_s$ the identity matrix of size $s$.

Blocking is a {\em batching} technique that can reduce the number of calls to the operator $A$ applied to individual column vectors,  maximize computational intensity by filling up the local cache with BLAS3 operations, and reduce the total number of sync points by performing inner products and normalization en masse \cite{BakDJ06, Bir15}.  In the context of Krylov subspaces, blocking can also lead to enriched subspaces by sharing information across column vectors instead of treating each right-hand side as an isolated problem.  How much information is shared across columns depends on the choice of block inner product.

Let $\spS$ be a $^*$-subalgebra of $\spCss$ with identity; i.e., $I \in \spS$ and when $ S,T \in \spS$, $\alpha \in \spC$, then $\alpha S +T, ST, S^* \in \spS$.
\begin{definition} \label{def:block_IPS}
	A mapping $\IPSnoarg$ from $\spC^{n\times s} \times \spCns$ to $\spS$ is called a {\em block inner product onto $\spS$} if it satisfies the following conditions for all $\vX,\vY,\vZ \in \spCns$ and $C \in \spS$:
	\begin{enumerate}[(i)]
		\item \textit{$\spS$-linearity}: $\IPS{\vX+\vY}{\vZ C} = \IPS{\vX}{\vZ} C + \IPS{\vY}{\vZ} C$; \label{S:linearity}
		\item \textit{symmetry}: $\IPS{\vX}{\vY} = \IPS{\vY}{\vX}^*$; \label{S:symmetry}
		\item \textit{definiteness}: $\IPS{\vX}{\vX}$ is positive definite if $\vX$ has full rank, and $\IPS{\vX}{\vX} = 0$ if and only if $\vX = 0$. \label{S:definiteness}
	\end{enumerate}
\end{definition}

\begin{definition} \label{def:scal_quot}
	A mapping $\N$ which maps all $\vX \in \spCns$ with full rank on a matrix $\N(\vX) \in \spS$ is called a {\em scaling quotient} if for all such $\vX$, there exists $\vY \in \spCns$ such that $\vX = \vY \N(\vX)$ and $\IPS{\vY}{\vY} = I_s$. 
\end{definition}

The scaling quotient is closely related to the intraorthogonalization routine discussed in Section~\ref{sec:bgs}.  Block notions of orthogonality and normalization arise organically from Definitions~\ref{def:block_IPS} and \ref{def:scal_quot}.
\begin{definition} \label{def:block_ortho_norm}
	Let $\vX, \vY \in \spCns$ and $\{\vX_j \}_{j=1}^m \subset \spCns$.
	\begin{enumerate}[(i)]
		\item $\vX, \vY$ are {\em block orthogonal}, if $\IPS{\vX}{\vY} = 0_s$.
		\item $\vX$ is {\em block normalized} if $\N(\vX) = I_s$.
		\item $\vX_1, \ldots, \vX_m$ are {\em block orthonormal} if $\IPS{\vX_i}{\vX_j} = \delta_{ij} I_s$.
	\end{enumerate}
\end{definition}

A set of vectors $\{\vX_j\}_{j=1}^m \subset \spCns$ {\em block spans} a space $\spK \subseteq \spCns$, and we write $\spK = \spnS\{\vX_j\}_{j=1}^m$ if
\[
\spK = \Big\{\sum_{j=1}^{m} \vX_j \Gamma_j : \Gamma_j \in \spS \mbox{ for } j = 1, \ldots, m \Big\}.
\]
The set $\{\vX_j\}_{j=1}^m$ constitutes a {\em block orthonormal basis} for $\spK = \spnS\{\vX_j\}_{j=1}^m$ if it is block orthonormal.

In this work, we consider only the {\em classical} and {\em global} block paradigms, described in Table~\ref{tab:block_inner_products}.  These paradigms represent the two extremes of information-sharing, with the classical approach maximizing information shared among columns and the global approach minimizing it; see, e.g., \cite[Theorem~3.3]{FroLS20}.  Moreover, the global paradigm leads to a lower complexity per iteration in Krylov subspace methods, because what are matrix-matrix products in the classical paradigm get reduced to scaling operations in the global one.  Many other paradigms are also possible; see, e.g., \cite{Dre20, DreEn19}.

\begin{table}[htbp!]
	\begin{center}
		\begin{tabular}{M{.2\textwidth}| M{.1\textwidth}| M{.25\textwidth}| M{.3\textwidth}}
			& $\spS$ &	$\IPS{\vX}{\vY}$ & $N(\vX)$ \\ \hline
			classical (\cl)	& $\spCss$	& $\vX^*\vY$	& $R$, where $\vX = \vQ R$, and $\vQ \in \spCns, \vQ^*\vQ = I_s$ \\ \hline
			global	(\gl) & $\spC I_s$	& $\tfrac{1}{s}\trace{(\vX^*\vY)}I_s$	& $\tfrac{1}{\sqrt{s}}\normF{\vX}I_s$
		\end{tabular}
	\end{center}
	\caption{\updated{Choices of $\spS$, $\IPSnoarg$, and $N$ in the classical and global block paradigms.} \label{tab:block_inner_products}}
\end{table}

\subsection{Block Gram-Schmidt} \label{sec:bgs}
Block Gram-Schmidt (BGS) is a routine for orthonormalizing a set of block vectors $\{\vX_j \}_{j=1}^m \subset \spCns$.  Writing
\[
\bXX := \begin{bmatrix} \vX_1 & \vX_2 & \cdots & \vX_m\end{bmatrix} \in \spC^{n \times ms},
\]
we define a BGS method as one that returns a block orthonormal $\bQQ \in \spC^{n \times ms}$ and a block upper triangular $\RR \in \spC^{ms \times ms}$ such that $\bXX = \bQQ \RR$.  Important measures in the analysis of BGS methods are the condition number of $\bXX$,
\begin{equation} \label{eq:cond(X))}
	\kappa(\bXX) := \frac{\sigma_{\max}(\bXX)}{\sigma_{\min}(\bXX)},
\end{equation}
i.e., the ratio between the largest and smallest singular values of $\bXX$, and the {\em loss of orthogonality (LOO)},
\begin{equation} \label{eq:loo}
	\norm{I - \IPS{\bQQ}{\bQQ}},
\end{equation}
where $\IPSnoarg$ is a generalized inner product as described in Section~\ref{sec:block_krylov}.

When we discuss the stabiliy of BGS methods, we refer to bounds on the loss of orthogonality in terms of machine precision, $\eps$.  We assume IEEE double precision here, so $\eps = \bigO{10^{-16}}$.

For categorizing BGS variants, we recycle the skeleton-muscle notation from \cite{CarLRetal22, Hoe10}, where {\em skeleton} refers to the \textit{inter}orthogonalization routine between block vectors, and the {\em muscle} refers to the \textit{intra}orthogonalization routine between the columns of a single block vector.  As a prototype, consider the Block Modified Gram-Schmidt (\BMGS) skeleton, given by Algorithm~\ref{alg:BMGS}.  Here, \IOnoarg denotes a generic muscle that takes $\vX \in \spC^{n \times s}$ and returns $\vQ \in \spC^{n \times s}$ and $R \in \spC^{s \times s}$ such that $\IPS{\vQ}{\vQ} = I_s$ and $\vX = \vQ R$.  For the classical paradigm, this could be any \updated{implementation of a} QR factorization: a column-wise Gram-Schmidt routine, Householder QR (\HouseQR), Cholesky QR (\CholQR), etc.  As for the global paradigm, there is only one possible muscle, given by the global scaling quotient, which effectively reduces to normalizing block vectors with a scaled Frobenius norm.  \updated{Consequently, intraorthogonalization does not actually occur in the global paradigm, as the columns of block vectors are not orthogonalized with respect to one another at all.}
\begin{algorithm}[htbp]
	\caption{$[\bQQ, \RR] = \BMGS(\bXX)$ \label{alg:BMGS}}
	\begin{algorithmic}[1]
		\State{$[\vQ_1, R_{11}] = \IO{\vX_1}$}
		\For{$k = 1, \ldots, p-1$}
		\State{$\vW = \vX_{k+1}$}
		\For{$j = 1, \ldots, k$}
		\State{$R_{j,k+1} = \IPS{\vQ_j}{\vW}$}
		\State{$\vW = \vW - \vQ_j R_{j,k+1}$}
		\EndFor
		\State{$[\vQ_{k+1}, R_{k+1,k+1}] = \IO{\vW}$}
		\EndFor
		\State \Return{$\bQQ = [\vQ_1, \ldots, \vQ_p]$, $\RR = (R_{jk})$}
	\end{algorithmic}
\end{algorithm}

We regard a single call to either $\IPSnoarg$ or \IOnoarg as one sync point, which is only possible in practice if single-reduce algorithms like \CholQR \cite{YamNYetal15} or \TSQR / \AllReduceQR \cite{DemGHetal12, MorYZ12} are employed for \IOnoarg.

\subsection{Block Krylov subspace methods} \label{sec:block_krylov}
The {\em $m$th block Krylov subspace for $A$ and $\vB$ (with respect to $\spS$)} is defined as
\begin{equation} \label{eq:block_krylov}
	\spKS_m(A, \vB) := \spnS\{\vB, A\vB, \ldots, A^{m-1}\vB\}.
\end{equation}

Block Arnoldi is often used to compute a basis for $\spKS_m(A, \vB)$, and it is typically implemented with \BMGS as the skeleton; see Algorithm~\ref{alg:BMGS-Arnoldi}. \texttt{BMGS-Arnoldi} accrues a high number of sync points due to the inner \texttt{for}-loop, where an increasing number of inner products is performed per block column.

\begin{algorithm}[htbp]
	\caption{$[\bVV_{m+1}, \HH_{m+1,m}, \Beta] = \BMGS\texttt{-Arnoldi}(A, \vB, m)$ \label{alg:BMGS-Arnoldi}}
	\begin{algorithmic}[1]
		\State{$[\vV_1, \Beta] = \IO{\vB}$}
		\For{$k = 1, \ldots, m$}
		\State{$\vW = A \vV_k$}
		\For{$j = 1, \ldots, k$}
		\State{$H_{j,k} = \IPS{\vV_j}{\vW}$}
		\State{$\vW = \vW - \vV_j H_{j,k}$}
		\EndFor
		\State{$[\vV_{k+1}, H_{k+1,k}] = \IO{\vW}$}  \label{line:BMGS}
		\EndFor
		\State \Return{$\bVV_{m+1} = [\vV_1, \ldots, \vV_{m+1}]$, $\HH_{m+1,m} = (H_{jk})$, $\Beta$}
	\end{algorithmic}
\end{algorithm}

Performing $m$ steps of a block Arnoldi routine returns the {\em block Arnoldi relation}
\begin{equation} \label{eq:block_arnoldi}
	A \bVV_m = \bVV_m \HH_m + \vV_{m+1} H_{m+1,m},
\end{equation}
where $\bVV_m$ $\spS$-spans $\spKS_m(A, \vB)$ and $\HH_m$ denotes the $ms \times ms$ \updated{principal} submatrix of $\HH_{m+1,m}$.

\subsubsection{Block full orthogonalization methods with low-rank modifications} \label{sec:bfom_mod}
We define
\begin{equation} \label{eq:Xm}
	\vX_m := \bVV_m \big( \HH_m + \MM \big)^\inv \vEhat_1 \Beta,
\end{equation}
where $\vEhat_1 = \vehat_1 \otimes I_s$ is a standard block unit vector, as the {\em (modified) block full orthogonalization method (BFOM)} for approximating \eqref{eq:AX=B}.  When $\MM = 0$, we recover BFOM, which minimizes the error in the $A$-weighted $\spS$-norm for $A$ hermitian positive definite \cite{FroLS17}.  There are infinitely many choices for $\MM$, but perhaps only a few useful ones, some of which are discussed in \cite{FroLS20}.  We will concern ourselves here with just $\MM = \HH_m^\cinv \big(\vEhat_m H_{m+1,m}^* H_{m+1,m}\big) \vEhat_m^*$, which gives rise to a block generalized minimal residual method (BGMRES) \cite{Sim96, SimGa96a, SimGa96b}.  As in \cite{FroLS20}, we implement \updated{BGMRES} as a modified BFOM here, with an eye towards downstream applications like $f(A)\vB$ where the BFOM form is explicitly needed.  In practice, there may be computational savings with a \updated{less} modular implementation; see, e.g., \cite{Gutknecht2007, GutSc08, GutSc09}.

\subsubsection{Static restarting and cospatial factors} \label{sec:restarting}
Restarting is a well established technique for reconciling a growing basis with memory limitations.  Define the residual of \eqref{eq:Xm} as
\begin{equation} \label{eq:Rm}
	\vR_m := \vB - A \vX_m.
\end{equation}
The basic idea of restarts is to use $\vR_m$ to build a new Krylov subspace, which we then use to approximate the error $\vE_m := A^\inv \vB - \vX_m$, which solves $A \vE = \vR_m$ in exact arithmetic.  Building a new Krylov subspace from $\vR_m$ directly is not a great idea, because it would require an extra computation with $A$.  Furthermore, we need a cheap, accurate, and ideally locally computable way to approximate $\norm{\vR_m}$ from one cycle to the next in order to monitor convergence.  In \cite{FroLS20} a static restarting method for low-rank modified BFOM is introduced that satisfies these requirements.  By ``static," we mean the basis size $m$ is fixed from one restart cycle to the next, in contrast to adaptive or dynamic restart cycle lengths.  We restate \cite[Theorem~4.1]{FroLS20}, which enables an efficient residual approximation and restarting procedure.
\begin{theorem} \label{thm:restarts}
	Suppose $\MM = \vM \vEhat_m^*$, where $\vM \in \spC^{ms \times s}$ and $\vEhat_m = \vehat_m \otimes I_s$. Define $\vU_m := \bVV_{m+1} \begin{bmatrix} \vM \\ -H_{m+1,m} \end{bmatrix}$ and let $\vXi_m := (\HH_m + \MM)^\inv \vEhat_1 \Beta$ be the block coefficient vector for the approximation $\vX_m = \bVV_m \vXi_m$ \eqref{eq:Xm} of the system \eqref{eq:AX=B}.  With $\vR_m$ as in \eqref{eq:Rm} it then holds that
	\begin{equation} \label{eq:cosptial_UXi}
		\vR_m = \vU_m \big(\vEhat_m^* \vXi_m\big).
	\end{equation}
\end{theorem}
We refer to the $s \times s$ matrix $\vEhat_m^* \vXi_m$ as a {\em cospatial factor}, and \eqref{eq:cosptial_UXi} as the {\em cospatial residual relation}.  The term {\em cospatial} refers to the fact that the columns of $\vR_m$ and those of $\vU_m$ span the same space.  Moreover, in exact arithmetic, it is not hard to see that
\begin{equation} \label{eq:normF_Rm}
	\normF{\vR_m} = \normF{\begin{bmatrix} \vM \\ -H_{m+1,m} \end{bmatrix} \big(\vEhat_m^* \vXi_m\big)},
\end{equation}
and the right-hand term can be computed locally (and possibly redundantly on each processor) for $m \ll n$.

If the approximate residual norm does not meet the desired tolerance, then we can compute the Arnoldi relation for $\spK_m(A, \vU_m)$ to obtain $\bVV_{m+1}^{(2)}$, $\HH_m^{(2)}$, $H_{m+1,m}^{(2)}$, and $\Beta^{(2)}$, where the superscript here and later denotes association to the restarted Krylov subspace.  We then approximate $\vE_m$ as
\[
\vD_m := \bVV_m^{(2)} (\HH_m^{(2)} + \MM^{(2)})^\inv \vEhat_1 \Beta^{(2)} \big(\vEhat_m^* \vXi_m\big),
\]
and update $\vX_m$ as
\[
\vX_m^{(2)} := \vX_m + \vD_m.
\]
The process is repeated, applying Theorem~\ref{thm:restarts} iteratively, until the desired residual tolerance is reached.

\begin{remark} 
	The analysis in \cite{FroLS17, FroLS20} is carried out in exact arithmetic.  Therefore, when we replace Algorithm~\ref{alg:BMGS-Arnoldi} with low-sync versions in Section~\ref{sec:low_sync}, all the results summarized in this section still hold, because all block Gram-Schmidt variants generate the same QR factorization in exact arithmetic.
\end{remark}


\section{Low-synchronization variants of block Arnoldi} \label{sec:low_sync}
To distinguish between block Arnoldi variants, we default to the name of the underlying block Gram-Schmidt skeleton.  We specify a {\em configuration} as \texttt{ip}-$\texttt{skel}\circ(\texttt{musc})$: inner product, skeleton, and muscle, respectively.  This naturally leads to bit of an ``alphabet soup," for which we ask the reader's patience, as it is crucial to precisely define algorithmic configurations for benchmarking.  Please refer often to Table~\ref{tab:skeletons}, which summarizes acronyms for all the Gram-Schmidt skeletons we consider in this text.  Note that the coefficient in front of the number of sync points per cycle is often used to describe low-sync methods; e.g., \BCGSPIP is a ``one-sync" method, while \BMGSSVL is a ``three-sync" method.

\begin{table}[htbp!]
	\begin{center}
		\begin{tabular}{M{.15\textwidth}|M{.25\textwidth}|M{.1\textwidth}|M{.15\textwidth}|M{.175\textwidth}}
			\thead{Underlying\\Gram-Schmidt\\skeleton\vspace{-2ex}} & \thead{Meaning\\behind\\abbreviations}		& \thead{\\Section}	& \thead{number of\\sync points\\per \updated{$m$-}cycle} & \thead{\updated{bound on}\\loss of\\orthogonality,\\\updated{assumption on $\kappa$}} \\ \midrule
			\BMGS					& Block Modified Gram-Schmidt 						& \ref{sec:block_krylov}	& \updated{$\frac{m(m+1)}{2}$}	& $\bigO{\eps} \kappa$, $\bigO{\eps} \kappa < 1$ \\ \hline
			\BCGSPIP				& Block Classical GS, Pythagorean with Inner Product& \ref{sec:bcgs_p}			& \updated{$m+1$}				& $\bigO{\eps} \kappa^2$, $\bigO{\eps} \kappa^2 < 1$ \\ \hline
			\updated{\BCGSPIO}		& \updated{Block Classical GS, Pythagorean with Intraorthogonalization}& \updated{\ref{sec:bcgs_p}}	& \updated{$2m+1$}				& \updated{$\bigO{\eps} \kappa^2$, $\bigO{\eps} \kappa^2 < 1$} \\ \hline
			\BMGSSVL~/ \BMGSLTS		& Schreiber \& Van Loan~/ Lower Triangular Solve 	& \ref{sec:svl_lts}			& \updated{$3m$} 							& $\bigO{\eps} \kappa$, $\bigO{\eps} \kappa < 1$ \\ \hline
			\BMGSCWY~/ \BMGSICWY	& Compact WY~/ Inverse Compact WY 					& \ref{sec:cwy_icwy}		& \updated{$m+2$} 							& \updated{--} \\ \hline
			\BCGSIROLS				& Inner Reorthogonalization (+), Low-Sync 			& \ref{sec:bcgs_iro_ls}		& \updated{$m+2$}							& \updated{--}
		\end{tabular}
	\end{center}
	\caption{Acronyms for BGS skeletons. \updated{Here ``$m$-cycle"} refers to a restart cycle, or the construction of \updated{$\bVV_{m+1}$}.  Loss of orthogonality is defined in \eqref{eq:loo}, and here $\kappa$ is shorthand for $\kappa([\vB \,\, A \bVV_m])$.  The loss of orthogonality \updated{bound for \BMGSLTS is conjecture and for \BMGSCWY, \BMGSICWY, and \BCGSIROLS, unknown}. \label{tab:skeletons}}
\end{table}

\begin{remark}
	The methods presented here are closely related to but not quite the same as the block methods used by Yamazaki et al.\ in \cite{YamTHetal20}, where \BMGS, \BCGSPIP, and \BCGSIROLS are employed as Gram-Schmidt skeletons in $s$-step Arnoldi (also known as communication-avoiding Arnoldi) \cite{BalCDetal14, Car15, Hoe10}, which is used to solve a linear system with a single right-hand side.  Recall that we are solving \eqref{eq:AX=B}, i.e., multiple right-hand sides simultaneously.
\end{remark}

\begin{remark} \label{rem:pseudocode}
	In the pseudocode for each algorithm, intermediate quantities like $\vW$ and $\vU$ are defined explicitly each iteration for readability.  In general, we purposefully avoid redefining quantities in a given iteration and instead only set an output (i.e., entries in $\Beta$, $\bVV_m$, or $\HH_{m+1,m}$) once all computations pertaining to that value are complete.  This approach simplifies mathematical analysis.  Exceptions include Algorithms~\ref{alg:BMGS} and \ref{alg:BMGS-Arnoldi}, where $\vW$ is redefined inside the \texttt{for}-loop as projected components are subtracted away from it.  In practice, it is preferable to save storage by overwriting block vectors of $\bVV_m$ instead of allocating separate memory for $\vW$ and $\vU$, for which there anyway may not be space.
\end{remark}

\subsection{\BCGSPIP and \BCGSPIO} \label{sec:bcgs_p} 
A simple idea for reducing the number of sync points in \BMGS is to condense the \texttt{for}-loop in lines~4-7 of Algorithm~\ref{alg:BMGS-Arnoldi} into a single inner product and subtraction,
\begin{align*}
	\HH_{1:k,k} & = \IPS{\bVV_k}{\vW} \\
	\vW 		& = \vW - \vW \HH_{1:k,k}
\end{align*}
This exchange gives rise to what is commonly referred to as the block classical Gram-Schmidt (\BCGS) method.  It is, however, rather unstable, with a loss of orthogonality worse than $\bigO{\eps} \kappa^2([\vB \,\, A\bVV_m])$ \cite{CarLR21}.  However, by making a correction based on the\updated{ block Pythagorean theorem (as derived in, e.g., \cite[Section~2.1]{CarLR21})}, we can guarantee a loss of orthogonality bounded by $\bigO{\eps} \kappa^2([\vB \,\, A\bVV_m])$, as long as $\bigO{\sqrt{\eps}} \kappa([\vB \,\, A\bVV_m]) \leq 1$.

One version of the corrected algorithm is given as Algorithm~\ref{alg:BCGSPIP-Arnoldi}.  The acronym ``PIP" stands for ``Pythagorean (variant) with Inner Product," due to how the factor $H_{k+1,k}$ is computed.  An alternative formulation based off \BCGSPIO (where ``PIO" stands for ``Pythagoren with IntraOrthogonalization") is also possible and is given as Algorithm~\ref{alg:BCGSPIO-Arnoldi}.  Note that in line~5, we use $\sim$ to denote that a full block vector need not be computed or stored here, just the $2s \times 2s$ scaling quotient $\Omega$. For subtle reasons, \BCGSPIO appears to be less reliable in practice (see Section~\ref{sec:stab_repro}).

\begin{algorithm}[H]
	\caption{$[\bVV_{m+1}, \HH_{m+1,m}, \Beta] = \BCGSPIP\texttt{-Arnoldi}(A, \vB, m)$ \label{alg:BCGSPIP-Arnoldi}}
	\begin{algorithmic}[1]
		\State{$[\vV_1, \Beta] = \IO{\vB}$}
		\For{$k = 1, \ldots, m$}
		\State{$\vW = A \vV_k$}
		\State{$\begin{bmatrix}	\HH_{1:k,k} \\ \Omega \end{bmatrix}
			= \IPS{\begin{bmatrix} \bVV_k & \vW \end{bmatrix}}{\vW}$}
		\State{$H_{k+1,k} = \chol(\Omega - \HH_{1:k,k}^* \HH_{1:k,k})$}
		\State{$\vV_{k+1} = \big(\vW - \bVV_k \HH_{1:k,k}\big) H_{k+1,k}^\inv$}
		\EndFor
		\State \Return{$\bVV_{m+1} = [\vV_1, \ldots, \vV_{m+1}]$, $\HH_{m+1,m} = (H_{jk})$, $\Beta$}
	\end{algorithmic}
\end{algorithm}

\begin{algorithm}[H]
	\caption{$[\bVV_{m+1}, \HH_{m+1,m}, \Beta] = \BCGSPIO\texttt{-Arnoldi}(A, \vB, m)$ \label{alg:BCGSPIO-Arnoldi}}
	\begin{algorithmic}[1]
		\State{$[\vV_1, \Beta] = \IO{\vB}$}
		\For{$k = 1, \ldots, m$}
		\State{$\vW = A \vV_k$}
		\State{$\HH_{1:k,k} = \updated{\IPS{\bVV_k}{\vW}}$}
		\State{$\left[\sim, \begin{bmatrix} W & 0 \\ 0 & H	\end{bmatrix}\right] = \IO{\begin{bmatrix} \vW & 0 \\ 0 & \HH_{1:k,k}  \end{bmatrix}}$} 
		\State{$H_{k+1,k} = \chol(W^* W - H^* H)$}
		\State{$\vV_{k+1} = \big(\vW - \bVV_k \HH_{1:k,k}\big) H_{k+1,k}^\inv$}
		\EndFor
		\State \Return{$\bVV_{m+1} = [\vV_1, \ldots, \vV_{m+1}]$, $\HH_{m+1,m} = (H_{jk})$, $\Beta$}
	\end{algorithmic}
\end{algorithm}

\subsection{\BMGSSVL / \BMGSLTS} \label{sec:svl_lts} 
Barlow developed and analyzed one of the first stabilized low-sync Gram-Schmidt methods by using the Schreiber-Van Loan representation of products of Householder transformations \cite{Bar19, SchVa89}.  Under modest conditions, this method-- which we denote here as \BMGSSVL-- has loss of orthogonality like \BMGS.  \updated{Its success depends on tracking the loss of orthogonality via an auxiliary matrix $\TT$ (as defined in lines~1, 2, 6, and 9 of Algorithm~\ref{alg:BMGSSVL-Arnoldi}) and using this matrix to make corrections each iteration.  A closely related method is \BMGSLTS, which is identical to \BMGSSVL except that the $\TT$ matrix is formed via lower-triangular solves instead of matrix products.  A column version of \BMGSLTS was first developed by \`{S}wirydowicz et al.~\cite{SwiLAetal20} and generalized to blocks by Carson et al.~\cite{CarLRetal22}. Although \BMGSLTS appears to behave identically to \BMGSSVL in practice,} a formal analysis for the former remains open.  We present Arnoldi versions of \BMGSSVL and \BMGSLTS as , with different colors highlighting the small differences between the methods.  In both methods, the main inner product in line~4 is performed as in \BCGS.  Meanwhile $\TT$ acts as a kind of buffer, storing the loss of orthogonality per iteration, which is used in successive iterations to make small corrections to the computation in line~4.  Balabanov and Grigori use a similar technique to stabilize randomized sketches of inner products, where instead of explicitly computing and storing $\TT$, they solve least squares problems to compute $\HH_{1:k,k}$ \cite{BalGr21, BalGr22}.

\begin{algorithm}[htbp!]
	\caption{$[\bVV_{m+1}, \HH_{m+1,m}, \Beta] = \textcolor{red}{\BMGSSVL}/\textcolor{blue}{\BMGSLTS}\texttt{-Arnoldi}(A, \vB, m)$ \label{alg:BMGSSVL-Arnoldi}}
	\begin{algorithmic}[1]
		\State{\updated{$\TT = I_{ms}$}}
		\State{$[\vV_1, \Beta, T_{11}] = \IO{\vB}$}
		\For{$k = 1, \ldots, m$}
		\State{$\vW = A \vV_k$}
		\State{$\vY = \IPS{\bVV_k}{\vW}$}
		\State{$\HH_{1:k,k} = \textcolor{red}{\TT_{1:k,1:k}^* \vY}
			\texttt{~OR~}
			\textcolor{blue}{\TT_{1:k,1:k}^\cinv \vY}$}
		\State{$[\vV_{k+1}, H_{k+1,k}, T_{k+1,k+1}] = \IO{\vW - \bVV_k \HH_{1:k,k}}$}
		\State{$\vZ = \IPS{\bVV_k}{\vV_{k+1}}$}
		\State{$\TT_{1:k,k+1} = \textcolor{red}{-\TT_{1:k,1:k} \vZ \, T_{k+1,k+1}}
			\texttt{~OR~}
			\textcolor{blue}{\vZ \, T_{k+1,k+1}}$}
		\EndFor
		\State \Return{$\bVV_{m+1} = [\vV_1, \ldots, \vV_{m+1}]$, $\HH_{m+1,m} = (H_{jk})$, $\Beta$}
	\end{algorithmic}
\end{algorithm}

\subsection{\BMGSCWY / \BMGSICWY} \label{sec:cwy_icwy}
A column-wise version of this algorithm was first presented by \'{S}wirydowicz et al.\ as \cite[Algorithm~8]{SwiLAetal20}.  To the best of our knowledge, we are the first to develop a block-wise formulation, which we refer to here as \texttt{BMGS-CWY-Arnoldi}, where CWY stands for ``compact WY," an alternative way to represent Householder transformations used to originally derive this algorithm.  A related Arnoldi algorithm, not treated in either \cite{SwiLAetal20} or \cite{ThoCRetal22}, is based on the inverse CWY (ICWY) form, and is given simultaneously with \BMGSCWY in Algorithm~\ref{alg:BMGSCWY-Arnoldi}.

It is important to note that \texttt{BMGS-CWY-Arnoldi} would not reduce to \cite[Algorithm~8]{SwiLAetal20} or \cite[Algorithm~6.1]{ThoCRetal22} for $s=1$, as we have one total sync point, due to the lack of a reorthonormalization step for $\vV_k$.  Algorithm~\ref{alg:BMGSCWY-Arnoldi} was largely derived by transforming \BMGSCWY and \BMGSICWY from \cite{CarLRetal22} into a block Arnoldi routine.  The most challenging part is tracking how the $\RR$ factor in the Gram-Schmidt formulation maps to $\HH_{m+1,m}$ and determining where to scale by the off-diagonal entry $H_{k,k-1}$ each iteration.  It is also possible to compute only with $\RR$ and reconstruct $\HH_{m+1,m}$ after $\bVV_{m+1}$ is finished; this approach proved to be much less stable in practice, however, due to the growing condition number of $\RR$.
\begin{algorithm}[htbp!]
	\caption{$[\bVV_{m+1}, \HH_{m+1,m}, \Beta] = \textcolor{red}{\BMGSCWY}/\textcolor{blue}{\BMGSICWY}\texttt{-Arnoldi}(A, \vB, m)$ \label{alg:BMGSCWY-Arnoldi}}
	\begin{algorithmic}[1]
		\State{$\TT = I_{(m+1)s}$}
		\State{$[\vV_1, \Beta] = \IO{\vB}$}
		\State{$\vU = \vV_1$}
		\For{$k = 1, \ldots, m+1$}
		\State{$\vW = A \vU$}
		\If{$k = 1$}
		\State{$H_{1,1} = \IPS{\vU}{\vW}$}
		\State{$\vU = \vW - \vV_1 H_{1,1}$}
		\Else
		\State{$\begin{bmatrix} \vY & \vZ \\ \Omega & \widetilde{P}	\end{bmatrix}
			= \IPS{\begin{bmatrix} \bVV_{k-1} & \vU \end{bmatrix}}{\begin{bmatrix} \vU & \vW \end{bmatrix}}$}
		\State{$H_{k,k-1} = \chol(\Omega)$}
		\State{$P = H_{k,k-1}^\cinv \widetilde{P}$}
		\State{$\TT_{1:k-1,k} = \textcolor{red}{-\TT_{1:k-1,1:k-1} \big(\vY H_{k,k-1}^\inv\big)}
			\texttt{~OR~}
			\textcolor{blue}{\vY H_{k,k-1}^\inv}$}
		\State{$H_{1:k,k} = \textcolor{red}{\TT_{1:k,1:k}^* \left(\begin{bmatrix} \vZ \\ P \end{bmatrix} H_{k,k-1}^\inv \right)}
			\texttt{~OR~}
			\textcolor{blue}{\TT_{1:k,1:k}^\cinv \left(\begin{bmatrix} \vZ \\ P \end{bmatrix} H_{k,k-1}^\inv \right)}$}
		\EndIf
		\State{$\vV_k = \vU H_{k,k-1}^\inv$}
		\State{$\vU = \vW H_{k,k-1}^\inv - \bVV_{k-1} H_{1:k,k}$}
		\EndFor
		\State \Return{$\bVV_{m+1} = [\vV_1, \ldots, \vV_{m+1}]$, $\HH_{m+1,m} = (H_{jk})$, $\Beta$}
	\end{algorithmic}
\end{algorithm}

\subsection{\BCGSIROLS} \label{sec:bcgs_iro_ls} 
One of the most intriguing of all the low-sync algorithms is DCGS2 \cite{BieLTetal22}, referred to as \CGSIROLS in \cite{CarLRetal22}.  This algorithm is a reformulation of reorthogonalized CGS with a single sync point derived by ``delaying" normalization to the next iteration, where operations are batched in a kind of $s$-step approach (where $s = 2$).  The column-wise version exhibits $\bigO{\eps}$ loss of orthogonality; a rigorous proof of the backward stability bounds remains open, however.  The block version, \BCGSIROLS, does not exhibit perfect $\bigO{\eps}$ LOO; see numerical results in \cite{CarLRetal22}.

Bielich et al.\ present a column-wise Arnoldi based on DCGS2 as Algorithm~4 in \cite{BieLTetal22}.  Our Algorithm~\ref{alg:BCGSIROLS-Arnoldi} is a direct block generalization of this algorithm with slight reformulations to match the aesthetics of Algorithm~\ref{alg:BMGSCWY-Arnoldi} and principles of Remark~\ref{rem:pseudocode}.  Note that, as in Algorithm~\ref{alg:BMGSCWY-Arnoldi}, we are able to compute $\HH_m$ directly, but we must track an auxiliary matrix $\vJ$ and scale several quantities by $H_{k-1,k-2}$.  An alternative version of Algorithm~\ref{alg:BCGSIROLS-Arnoldi} based more directly on \BCGSIROLS from \cite[Algorithm~7]{CarLRetal22} is included in the code but not described here.

\begin{algorithm}[htbp!]
	\caption{$[\bVV_{m+1}, \HH_{m+1,m}, \Beta] = \BCGSIROLS\texttt{-Arnoldi}(A, \vB, m)$ \label{alg:BCGSIROLS-Arnoldi}}
	\begin{algorithmic}[1]
		\State{$[\vV_1, \Beta] = \IO{\vB}$}
		\State{$\vU = \vV_1$}
		\For{$k = 1, \ldots, m+1$}
		\State{$\vW = A \vU$}
		\If{$k = 1$}
		\State{$\vJ = \IPS{\vU}{\vW}$}
		\State{$H_{1,1} = \vJ$}
		\State{$\vU = \vW - \vV_1 \vJ$}
		\Else
		\State{$\begin{bmatrix} \vY & \vZ \\ \widetilde{\Omega} & \widetilde{P} \end{bmatrix}
			= \IPS{\begin{bmatrix} \bVV_{k-1} & \vU \end{bmatrix}}{\begin{bmatrix} \vU & \vW \end{bmatrix}}$}
		\State{$\Omega = \widetilde{\Omega} - \vY^* \vY$}
		\State{$H_{k,k-1} = \chol(\Omega)$}
		\State{$\HH_{1:k-1,k-1} = \vJ + \vY$}
		\State{$P = H_{k,k-1}^\cinv (\widetilde{P} - \vY^* \vZ)$}
		\State{$\vJ = \left(\begin{bmatrix} \vZ \\ P	\end{bmatrix} - \HH_{1:k,1:k-1} \vY \right) H_{k,k-1}^\inv$}
		\State{$\vV_{k} = (\vU - \bVV_{k-1} \vY) H_{k,k-1}^\inv$}
		\State{$\vU = \left(\vW - \bVV_k \begin{bmatrix} \vZ \\ P \end{bmatrix} \right) H_{k,k-1}^\inv$}
		\EndIf
		\EndFor
		\State \Return{$\bVV_{m+1} = [\vV_1, \ldots, \vV_{m+1}]$, $\HH_{m+1,m} = (H_{jk})$, $\Beta$}
	\end{algorithmic}
\end{algorithm}


\section{Adaptive restarting} \label{sec:stab_repro}
Reproducibility and stability are not mutually exclusive.  This realization is precisely the motivation for an adaptive restarting routine and can be demonstrated by a simple example.

Consider the \tridiag test case from Section~\ref{sec:tridiag} with $n = 100$.  Notably, both $A$ and $\vB$ are deterministic quantities; neither is defined with random elements.  In \MATLAB, it is possible to specify the number of threads on which a script is executed via the built-in \texttt{maxNumCompThreads} function.\footnote{\url{https://mathworks.com/help/matlab/ref/maxnumcompthreads.html}. Accessed 8 August 2022.}  We solve $A \vX = \vB$ with Algorithms~\ref{alg:BCGSPIP-Arnoldi} and~\ref{alg:BCGSPIO-Arnoldi} while varying the multithreading setting from 1 to 16 on a standard node of the Mechthild cluster; see the beginning of Section~\ref{sec:examples} for more details about the cluster.   For both algorithms, we employ a variant of \MATLAB's Cholesky routine \chol, which stores a flag when \chol determines a matrix is too ill-conditioned to be factorized.  This flag is fed to the linear solver driver of \LSBA (\texttt{bfom}), which halts the process when the flag is true.  Through the following discussion, we refer to this flag as the ``\nan-flag," \updated{because ignoring it leads to computations with ill-defined quantities.}

\begin{figure}[htbp!]
	\begin{tabular}{cc}
		\resizebox{.44\textwidth}{!}{\includegraphics{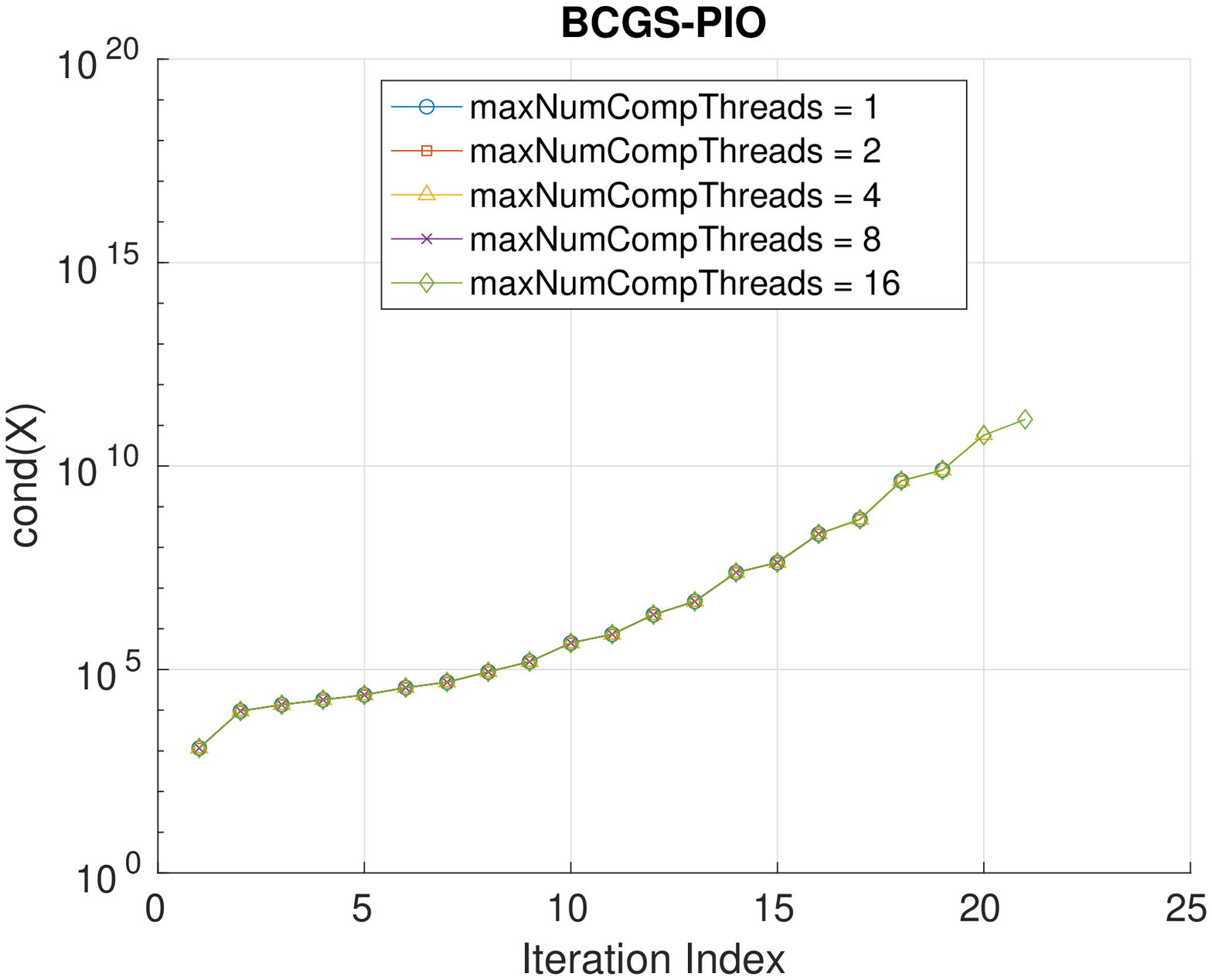}} &
		\resizebox{.44\textwidth}{!}{\includegraphics{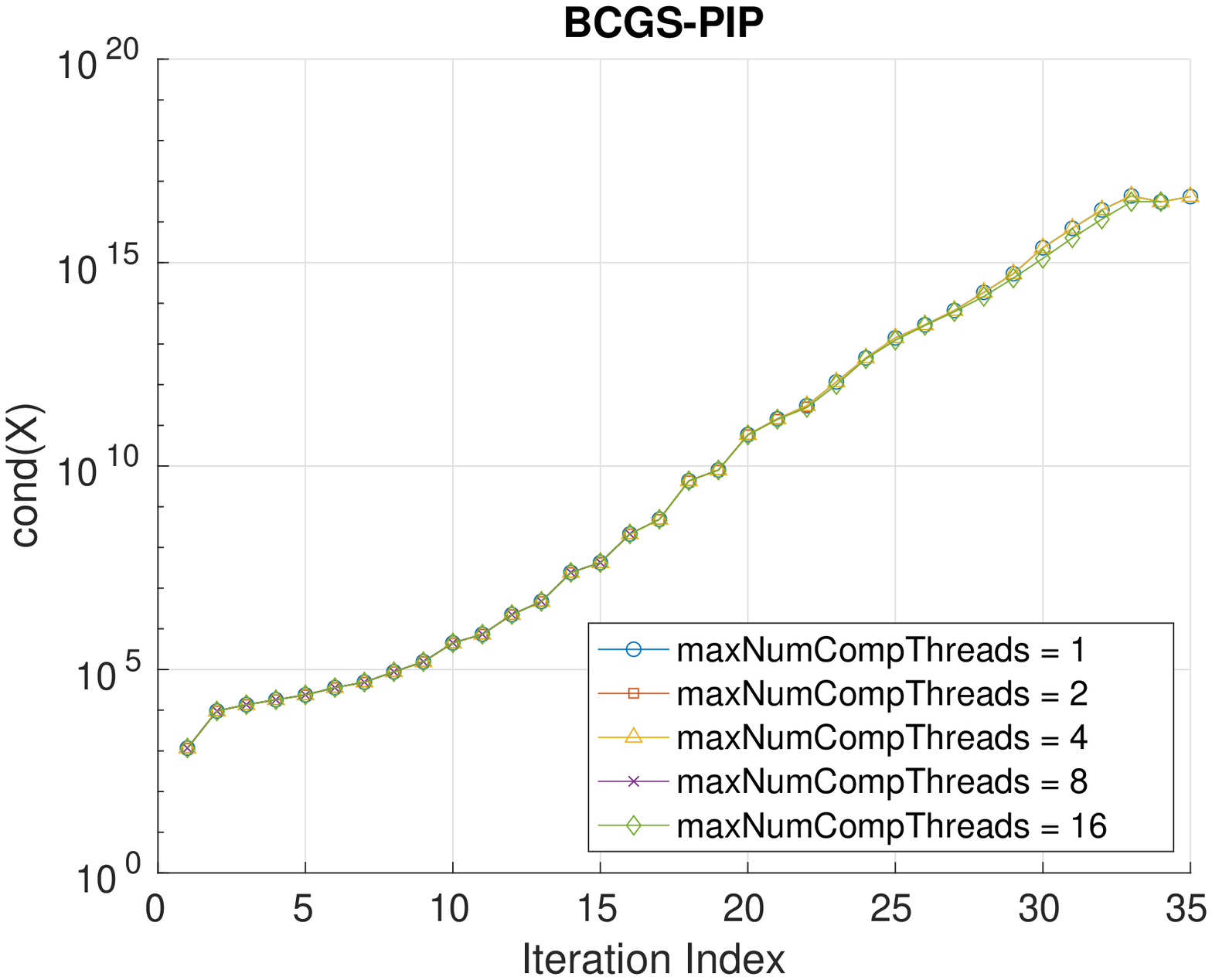}} \\
		& \\
		\resizebox{.44\textwidth}{!}{\includegraphics{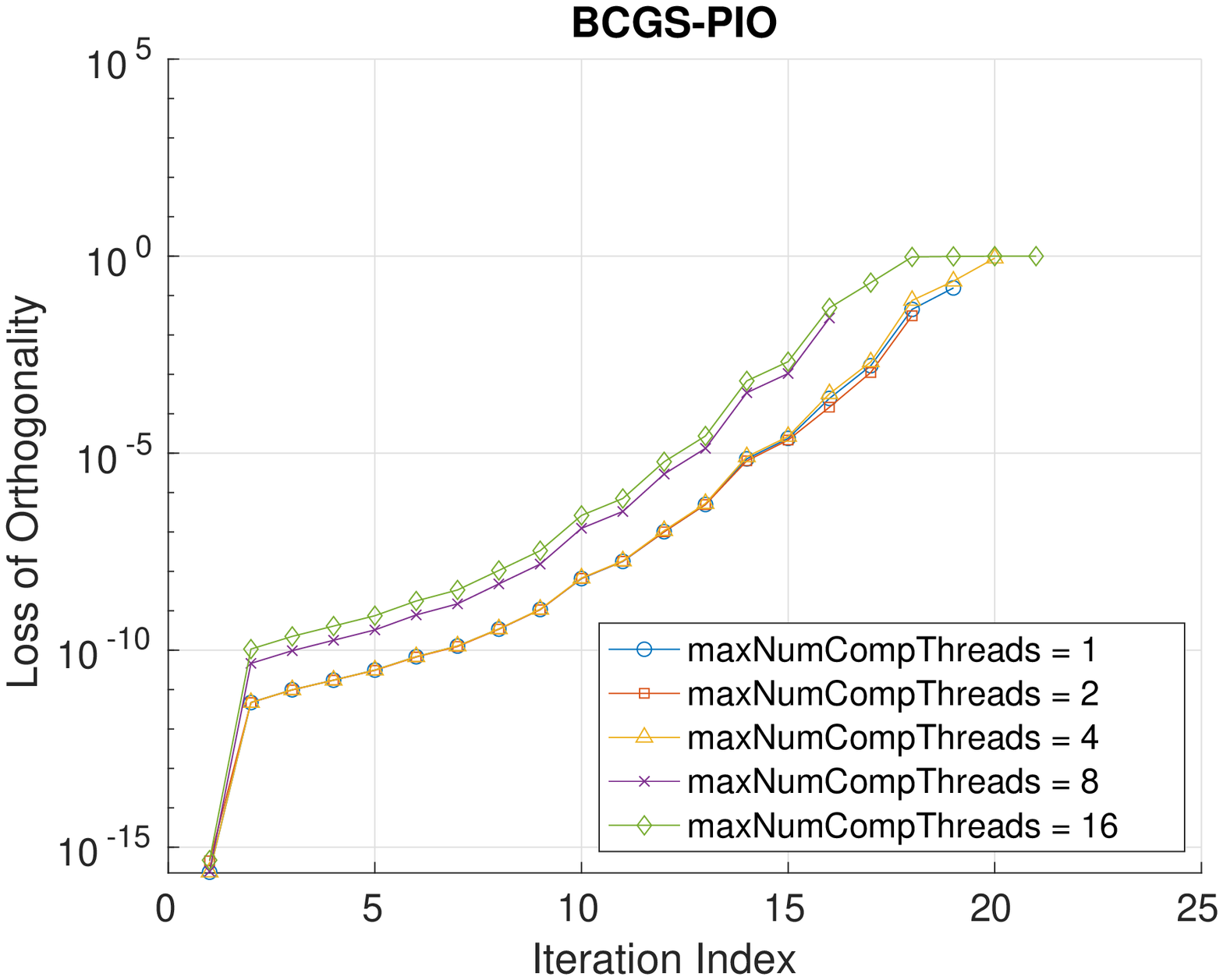}} &
		\resizebox{.44\textwidth}{!}{\includegraphics{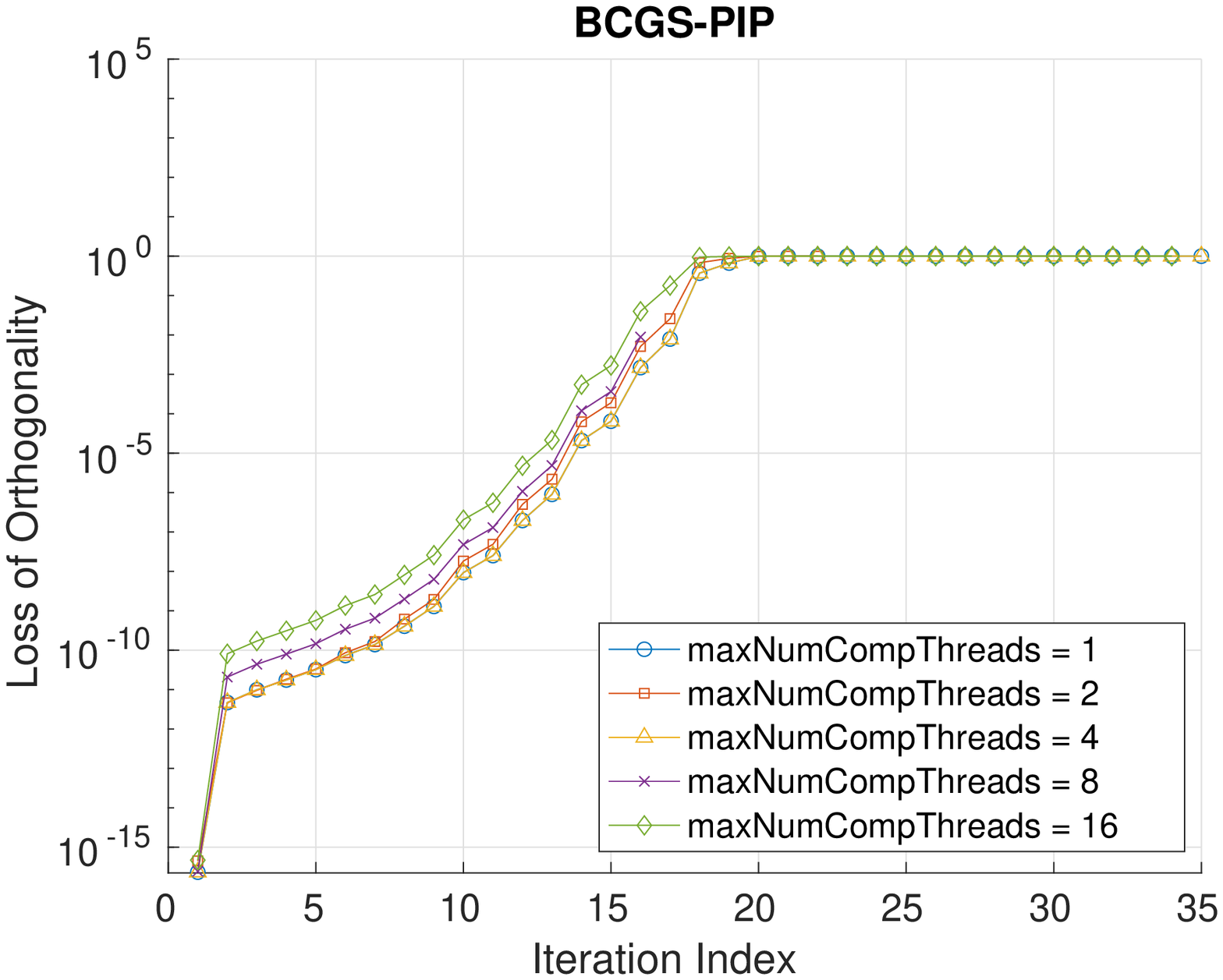}} \\	
	\end{tabular}
	\caption{Multithreading example for \tridiag problem with $n=100$, no restarts, and maximum basis size $m = 50$. \label{fig:tridiag_multithread}}
\end{figure}

Figure~\ref{fig:tridiag_multithread} displays the loss of orthogonality \updated{\eqref{eq:loo}} and $\kappa([\vB \,\, A \bVV_k])$ for different thread counts.  The condition numbers for all thread counts and both methods are hardly affected, except for some slight deviation for \BCGSPIP and 16 threads.  The LOO plots are more telling: for both methods, changing the thread count directly affects the LOO and how many iterations the method can compute before encountering a \nan-flag.  We allowed for a maximum basis size of $m=50$, but no method can compute that far.  \BCGSPIO with 8 threads gives up first at 16 iterations; \BCGSPIP with 1 and 4 threads makes it all the way to 35 iterations.  Among the \BCGSPIO methods, there are orders of magnitude differences between the attained LOO.

\updated{This situation is perplexing on the surface}: the problem is static, and the same code has been run every time.  The only variable is the thread count.

There are two subtle issues that affect reproducibility in this case: 1) the configuration of math kernel libraries according to the parameters of the operating system and hardware,\footnote{\url{https://www.intel.com/content/www/us/en/develop/documentation/onemkl-linux-developer-guide/top/obtaining-numerically-reproducible-results/reproducibility-conditions.html}. Accessed 8 August 2022.} and 2) guaranteed stability bounds.  As for stability bounds, it is important to note that both \BCGSPIO and \BCGSPIP have a complete backward stability analysis \cite{CarLR21}.  Both methods have $\bigO{\eps} \kappa^2([\vB \,\, A \bVV_k])$ loss of orthogonality, as long as $\kappa([\vB \,\, A \bVV_k]) \leq \bigO{\frac{1}{\sqrt{\eps}}} = \bigO{10^8}$ and as long as the \IOnoarg for \BCGSPIO behaves no worse than \CholQR.  (For this test, we used \HouseQR, \MATLAB's built-in \texttt{qr} routine, which is unconditionally stable and therefore behaves better than \CholQR \cite{Hig02}.)  For both methods, $\kappa([\vB \,\, A \bVV_k])$ exceeds $\bigO{10^8}$ around iteration 15.  At that point, the assumptions for the LOO bounds are no longer satisfied.  The fact that either algorithm continues to compute something useful after that point is a lucky accident.

Computing $\kappa([\vB \,\, A \bVV_k])$ every iteration to check whether the LOO bounds are satisfied is not practical.  We therefore propose a simple adaptive restarting regime based on whether \chol raises a \nan-flag, which happens whenever \chol \updated{is fed} a numerically non-positive definite matrix.  When a \nan-flag is raised, we give up computing a new basis vector and go back to the last safely computed basis vector, which is then used to restart.  Simultaneously, the maximum basis size $m$ is also reduced.  It is possible that an algorithm exhausts its maximum allowed restarts and basis size before converging; indeed, we have observed this often for \BCGSPIP in examples not reported here.  At the same time, there are many scenarios in which restarting is an adequate band-aid, thus allowing computationally cheap, one-sync algorithms line \BCGSPIP to salvage progress and converge, oftentimes faster than competitors.  See Section~\ref{sec:examples} for demonstrations.

\begin{remark}
	The restarted framework outlined in Section~\ref{sec:restarting} does not change fundamentally with adaptive cycle lengths; only the notation becomes more complicated.  We omit the details here.
\end{remark}

\section{Numerical benchmarks} \label{sec:examples}
Our treatment of BGS and block Krylov methods is hardly exhaustive.  It is not our goal to determine the optimal block Arnoldi configuration at this stage, but rather to demonstrate the functionality of a benchmarking tool for the fair comparison of possible configurations on different problems.  To this end, we restrict ourselves to the options below:
\begin{itemize}
	\item inner products: \cl\xspace (classical), \gl\xspace (global)
	\item skeletons: Table~\ref{tab:skeletons}
	\item muscles: \CholQR, which has $\bigO{\eps} \kappa^2$ loss of orthogonality guaranteed only for $\bigO{\eps} \kappa^2 < 1$, but is a simple, single-reduce algorithm.  In practice, we would recommend \TSQR/\AllReduceQR \cite{DemGHetal12, MorYZ12}, which has $\bigO{\eps}$ loss of orthogonality and the same number of sync points, but is difficult to program in \MATLAB due to limited parallelization and message-passing features.  Other low-sync muscles are programmed in \LSBA as well, and the user can easily integrate their own.  Note that \BCGSPIP does not require a muscle, and \BMGSCWY, \BMGSICWY, and \BCGSIROLS only call a muscle once, in the first iteration of a new basis. \BMGSSVL and \BMGSLTS are forced to use their column-wise counterparts \MGSSVL and \MGSLTS (both 3-sync), respectively, and global methods are forced to use the global muscle (i.e., \updated{normalization without intraorthogonalization via the scaled Frobenius norm}).
	\item modification: \texttt{none} (FOM), \texttt{harmonic} (GMRES)
\end{itemize}

All results are generated by the \LSBA \MATLAB package.  A single script (\texttt{paper\_script.m}) comprises all the calls for generating the results in this manuscript.  \LSBA is written as modularly as possible, to facilitate the exchange of inner products, skeletons, muscles, and modifications.  While the timings reported certainly do not reflect the optimal performance for any of the methods, they do reflect a fair comparison across implementations and provide insights for possible speed-ups when these methods are ported to more complex architectures.  The code is also written so that sync points (\innerprod and \intraortho) and other potentially communication-intensive operations (\matvec and \basiseval) are separate functions that can be tuned individually.

Every test script (including the example from Section~\ref{sec:stab_repro}) has been executed in \MATLAB R2019b on 16 threads of a single, standard node of Linux Cluster Mechthild at the Max Planck Institute for Dynamics of Complex Technical Systems in Magdeburg, Germany.\footnote{\url{https://www.mpi-magdeburg.mpg.de/cluster/mechthild}. Accessed 8 August 2022.}  A standard node comprises 2 Intel Xeon Silver 4110 (Skylake) CPUs with 8 Cores each (64KB L1 cache, 1024KB L2 cache), a clockrate of 2.1 GHz (3.0 GHz max), and 12MB shared L3 cache each.  We further focus on small problems that easily fit in the L3 Cache, which is easy to guarantee with sparse $A$, $n \leq 10^4$, and $s \leq 10$.  Given that the latency between CPUs on a single node is small relative to exascale machines, we expect small improvements observed in these test cases to translate to bigger gains in a more complex setting.

For the timings, we measure the total time spent to reach a specified error tolerance.  We run each test 5 times and average over the timings.  We also calculate several intermediate measures, namely counts for $A$-calls, applications of $\bVV_k$, and sync points.  In addition, we plot the convergence history in terms of the following quantities per iteration: relative residual, relative error, $\kappa([\vB \, \, A \bVV_k])$, and loss of orthogonality (LOO) \eqref{eq:loo}.  When a ground truth solution $\vX_*$ is provided, the error is calculated as
\[
\normF{\vX_k - \vX_*} / \normF{\vX_*},
\]
For all our examples, $\vX_*$ is computed by \MATLAB's built-in \texttt{backslash} operator.  The residual is approximated by \eqref{eq:normF_Rm} and is scaled by $\normF{\vB}$.  A summary of the parameters for all benchmarks can be found in Table~\ref{tab:parameters}.  Except for \tridiag and \lapl, all examples are taken from the SuiteSparse Matrix Collection \cite{DavH11}.  Via the \texttt{suite\_sparse.m} script, it is possible to run tests on any benchmark from this collection.
\begin{table}[htbp]
	\begin{center}
		\begin{tabular}{l|c|c|c|c|c|c}
			\thead{test name}	& $\kappa(A)$ 	 & \thead{$n$}	& \thead{$s$} & \thead{$m$}	& \thead{modification}	& \thead{tol} 	\\ \midrule
			\tridiag     		& $\bigO{10^3}$	 & $1000$		& $2$		  & $70$ 		& FOM				  	& $10^{-10}$		\\
			\powersystem 		& $\bigO{10^6}$	 & $1138$		& $5$		  & $30$		& GMRES		  			& $10^{-6}$		\\
			\texttt{circuit\_2} & $\bigO{10^5}$	 & $4510$ 		& $5$ 		  & $10$        & GMRES 			  	& $10^{-6}$		\\
			\texttt{rajat03}	& $\bigO{10^7}$	 & $7602$		& $5$		  & $10$		& GMRES				  	& $10^{-6}$		\\
			\texttt{Kaufhold}	& $\bigO{10^{14}}$ & $8765$		& $5$		  & $10$		& GMRES				  	& $10^{-6}$		\\
			\texttt{t2d\_q9}	& $\bigO{10^3}$	 & $9801$		& $5$		  & $10$		& GMRES				  	& $10^{-6}$		\\
			\lapl				& $\bigO{10^3}$	 & $10000$		& $10$		  & $25$		& FOM				  	& $10^{-6}$
		\end{tabular}
	\end{center}
	\caption{Test properties and parameter choices. \label{tab:parameters}}
\end{table}

\subsection{\tridiag} \label{sec:tridiag}
The operator $A$ is defined as a sparse, tridiagonal matrix with $1$ on the off-diagonals and $-1, -2, \ldots, -n$ on the diagonal, where $n$ is also the size of $A$.  Clearly $A$ is symmetric.  The right-hand side $\vB$ has two columns, where the first has identical elements $\frac{1}{\sqrt{n}}$ and the second is $1, 2, \ldots, n$.  This example is actually procedural, in the sense that a user can choose a desired $n$.  At the same time, a larger $n$ necessarily leads to a worse condition number.

\begin{figure}[htbp!]
	\resizebox{.95\textwidth}{!}{\includegraphics{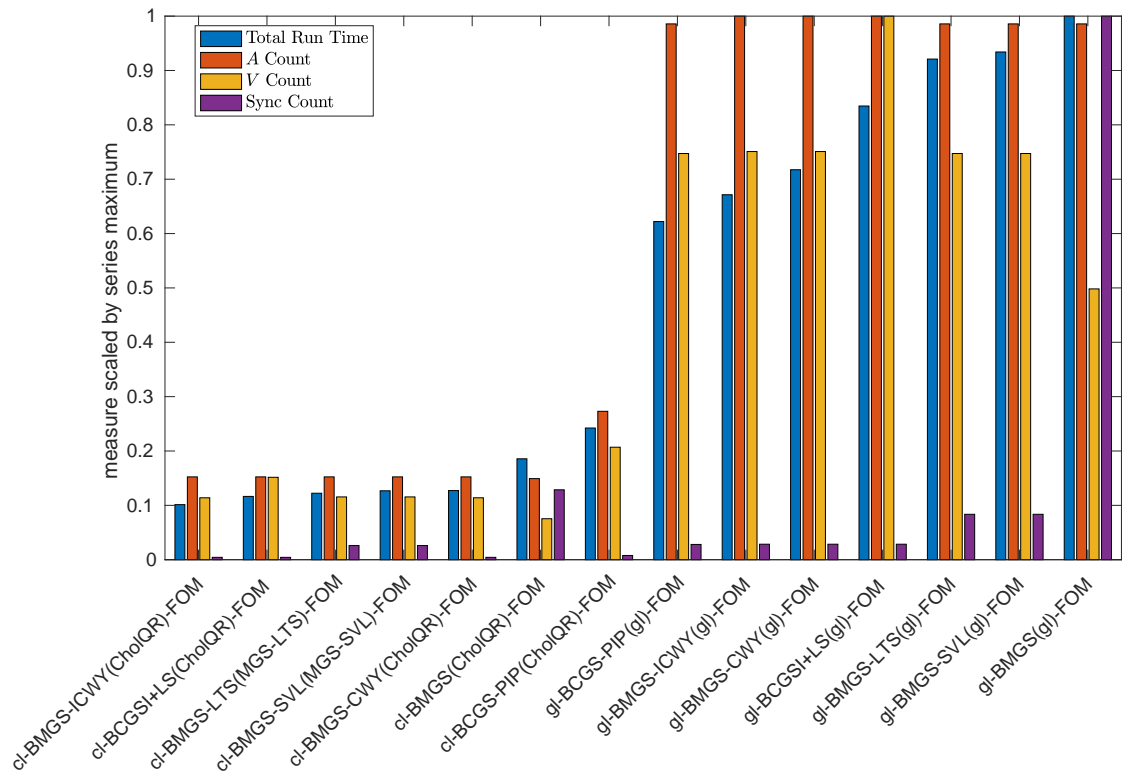}}
	\caption{Results from \tridiag example. \label{fig:tridiag}}
\end{figure}

Figure~\ref{fig:tridiag} presents the total run time per configuration as well as operator counts as a bar chart; see Table~\ref{tab:tridiag} in the Appendix for more details.  The fastest methods are the stabilized low-sync variants.  Despite being the computationally cheapest classical method per iteration, \cl-\BCGSPIP is notably slower than \cl-\BMGS, because its inherent instability requires restarting $3$ times (and therefore additional applications of $A$ and $\bVV_k$) before converging.  The method with the fewest $\bVV_k$ evaluations is \cl-\BMGS, which is to be expected, since the basis is split up and applied one block column at a time in the inner-most loop; see Algorithm~\ref{alg:BMGS}.

The fastest global method, \gl-\BCGSPIP, is significantly slower even than the slowest classical method.  In fact, all global methods require over 6 times as many total iterations as the fastest classical method to converge; this is in line with the theory of Section~\ref{sec:block_krylov}.  In this particular case, the floating-point savings per iteration do not outweigh the sheer amount of time needed for all the extra $A$-calls.  Nevertheless, the one-sync global methods (\gl-\BCGSPIP, \gl-\BMGSCWY, \gl-\BMGSICWY, and \gl-\BCGSIROLS) have relatively low sync counts, compared even to \cl-\BMGS.

\begin{figure}[htbp!]
	\begin{tabular}{cc}
		\resizebox{.44\textwidth}{!}{\includegraphics{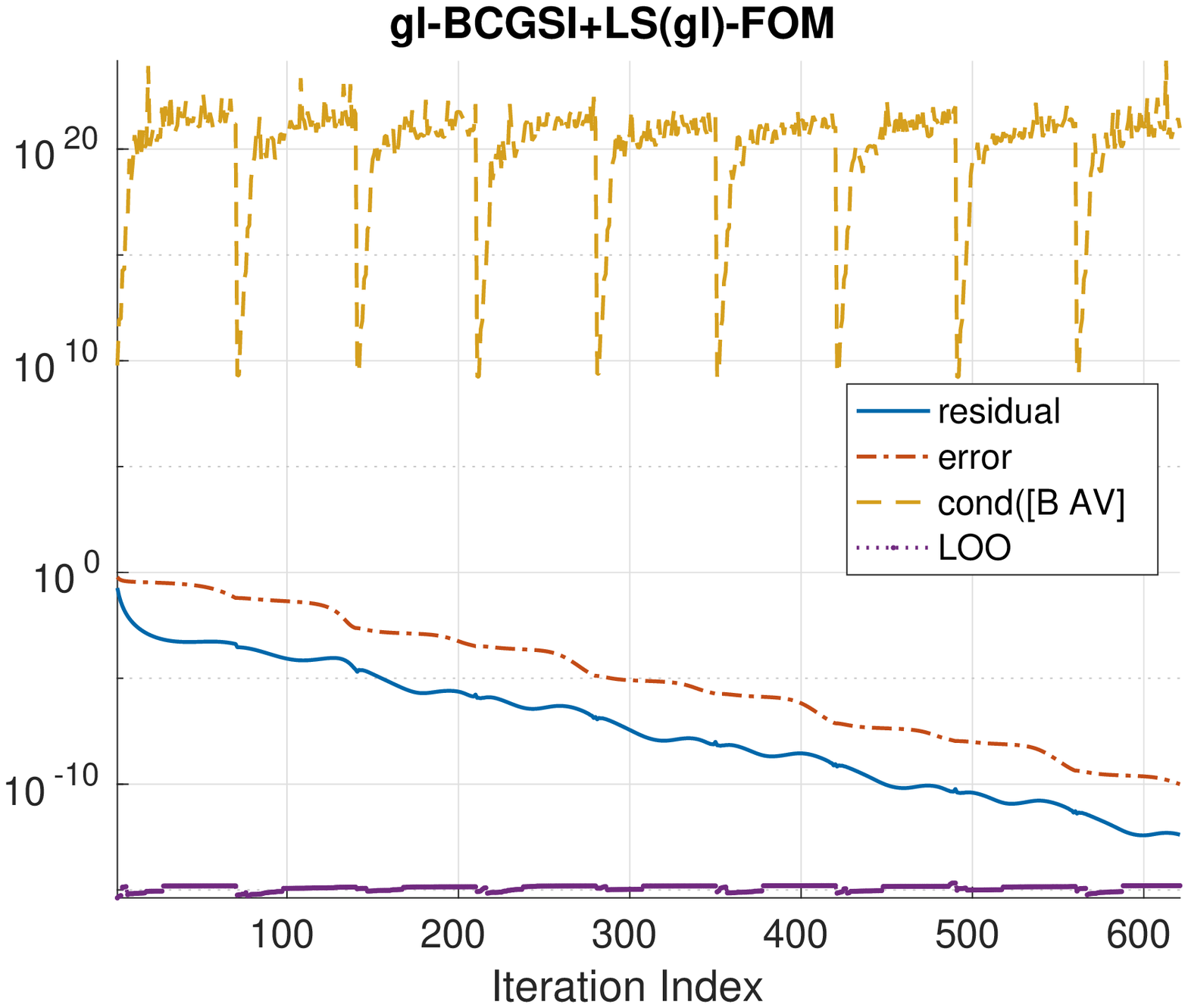}} &
		\resizebox{.44\textwidth}{!}{\includegraphics{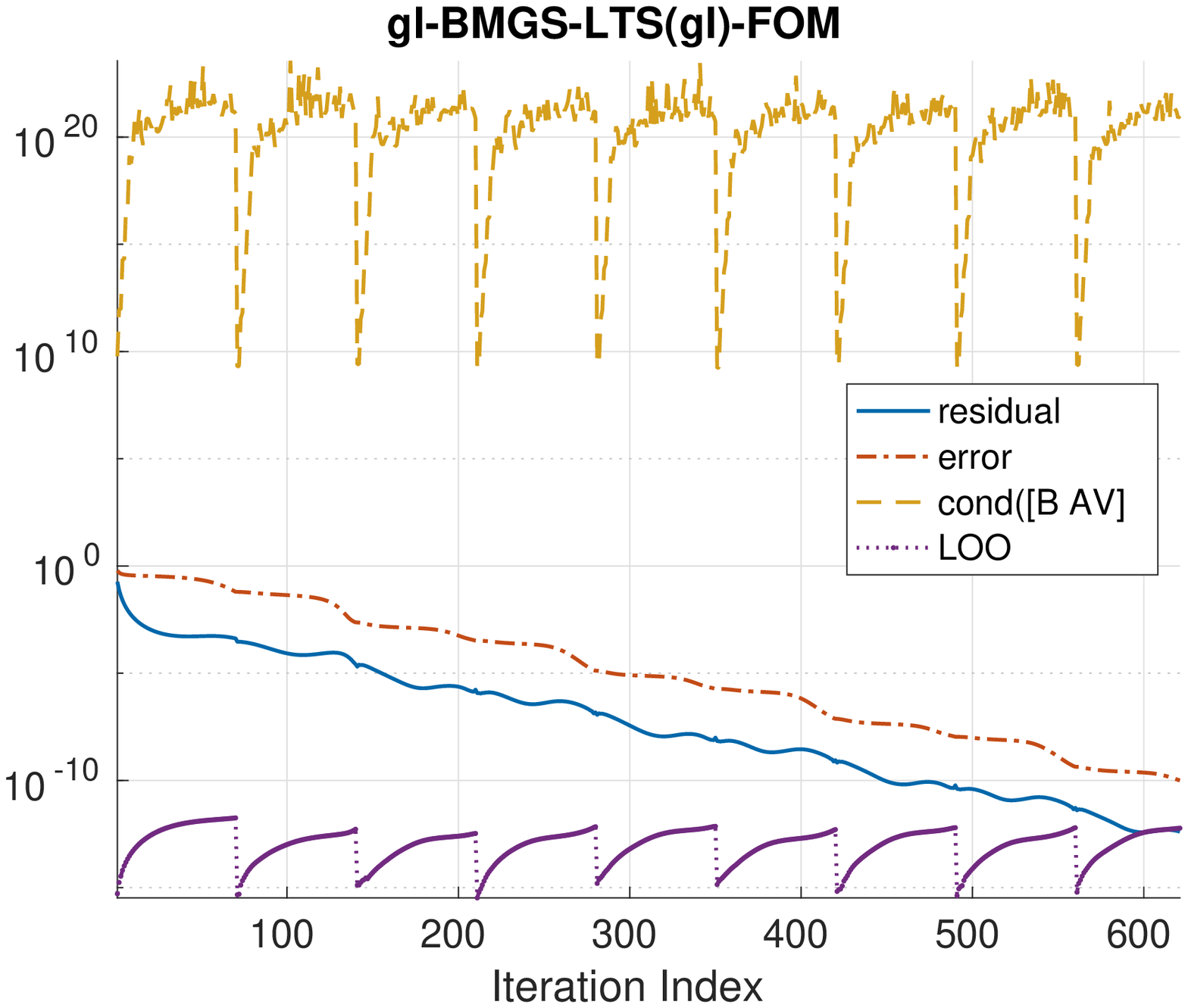}} \\
		& \\
		\resizebox{.44\textwidth}{!}{\includegraphics{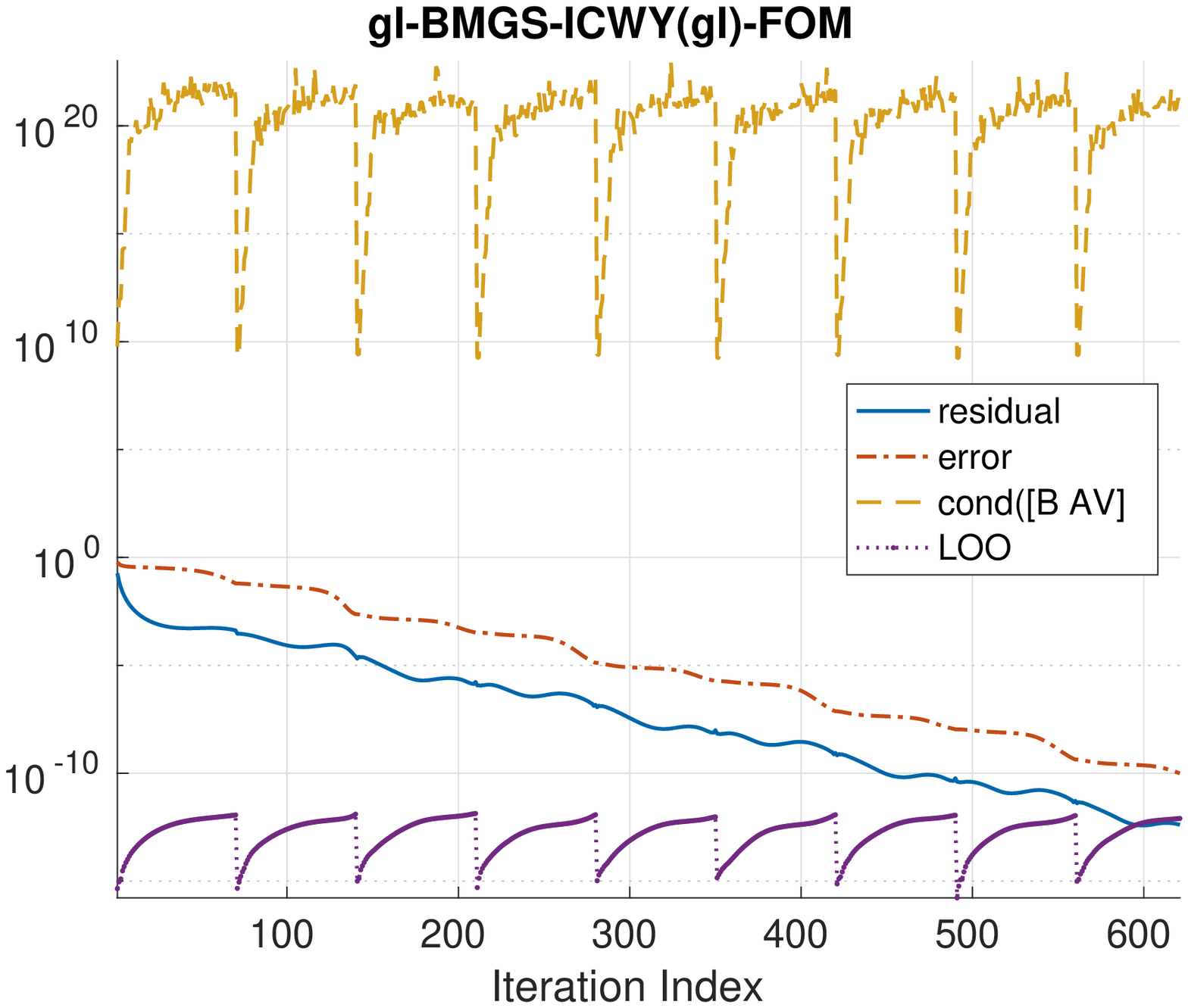}} &
		\resizebox{.44\textwidth}{!}{\includegraphics{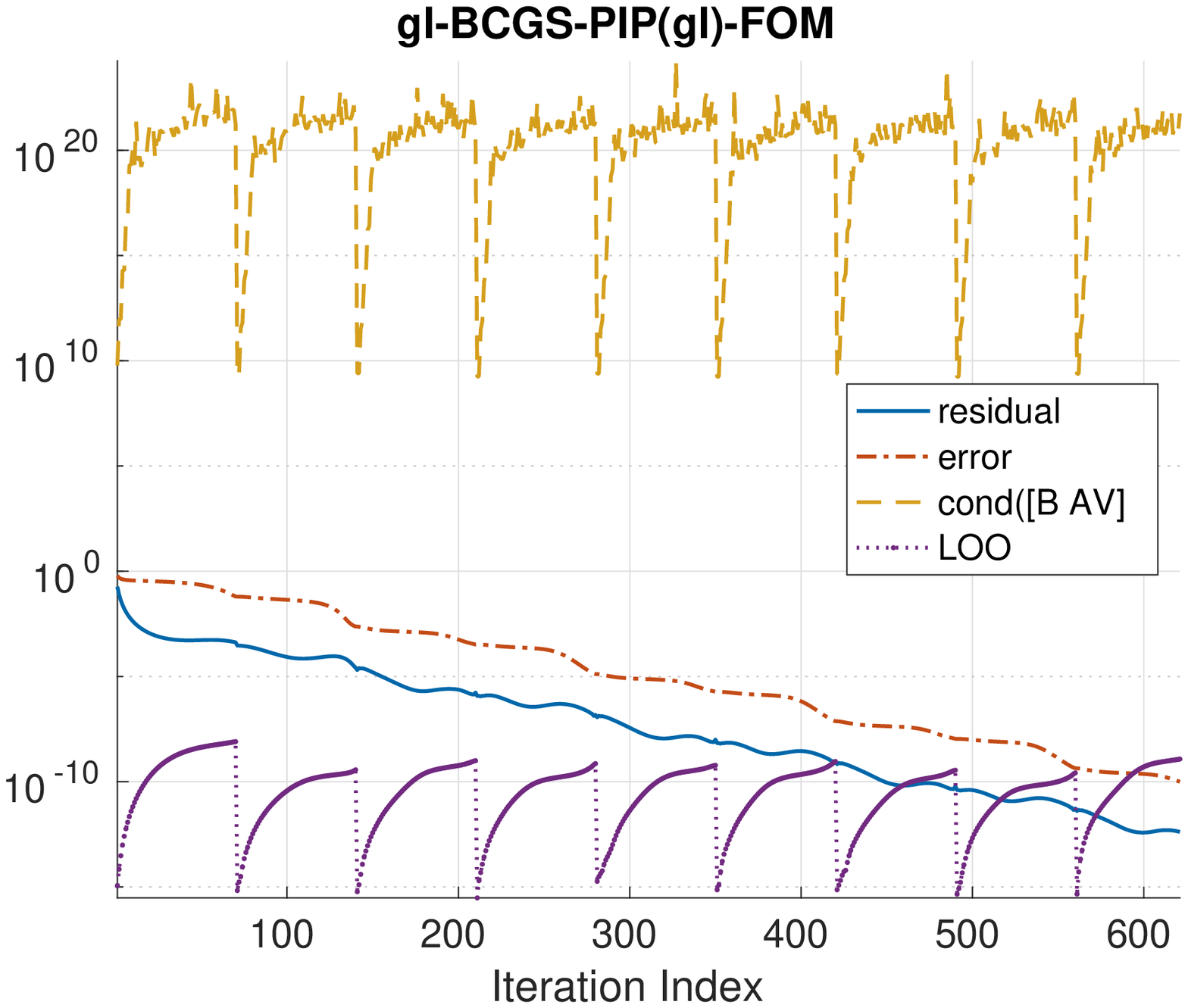}} \\
	\end{tabular}
	\caption{Convergence histories of some global variants for \tridiag example. \label{fig:tridiag_gl_conv}}
\end{figure}

\begin{figure}[htbp!]
	\begin{tabular}{cc}
		\resizebox{.44\textwidth}{!}{\includegraphics{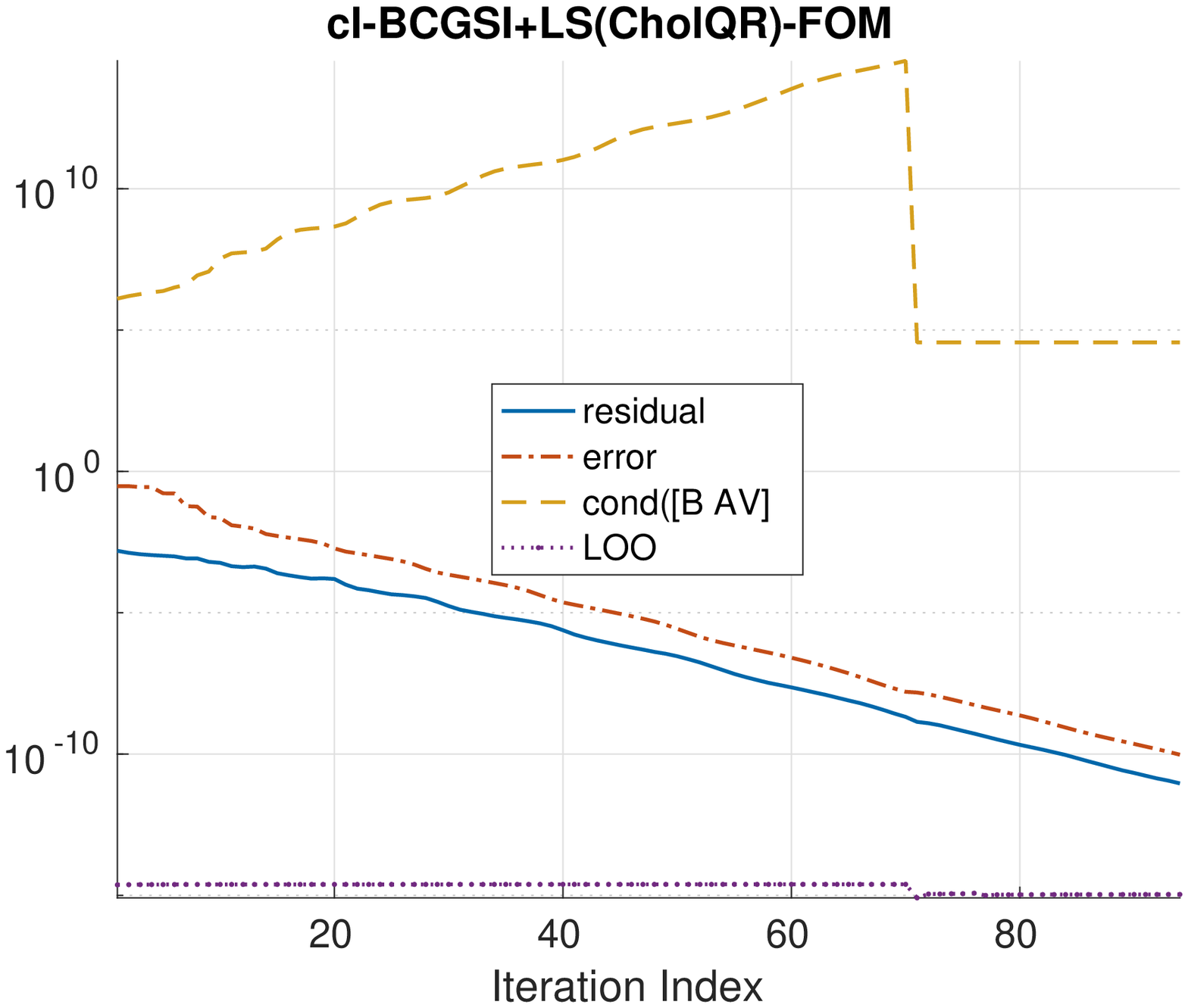}} &
		\resizebox{.44\textwidth}{!}{\includegraphics{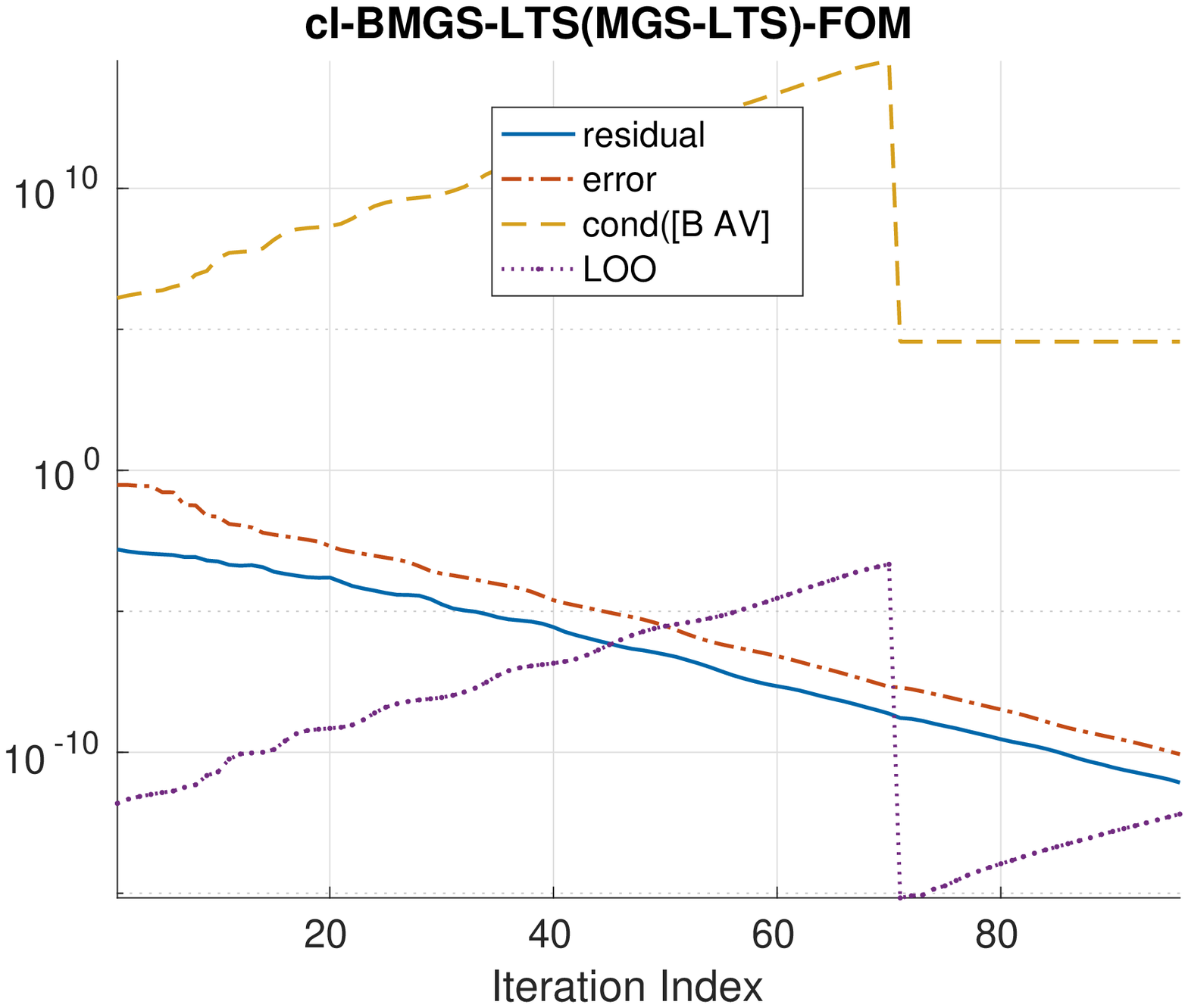}} \\
		& \\
		\resizebox{.44\textwidth}{!}{\includegraphics{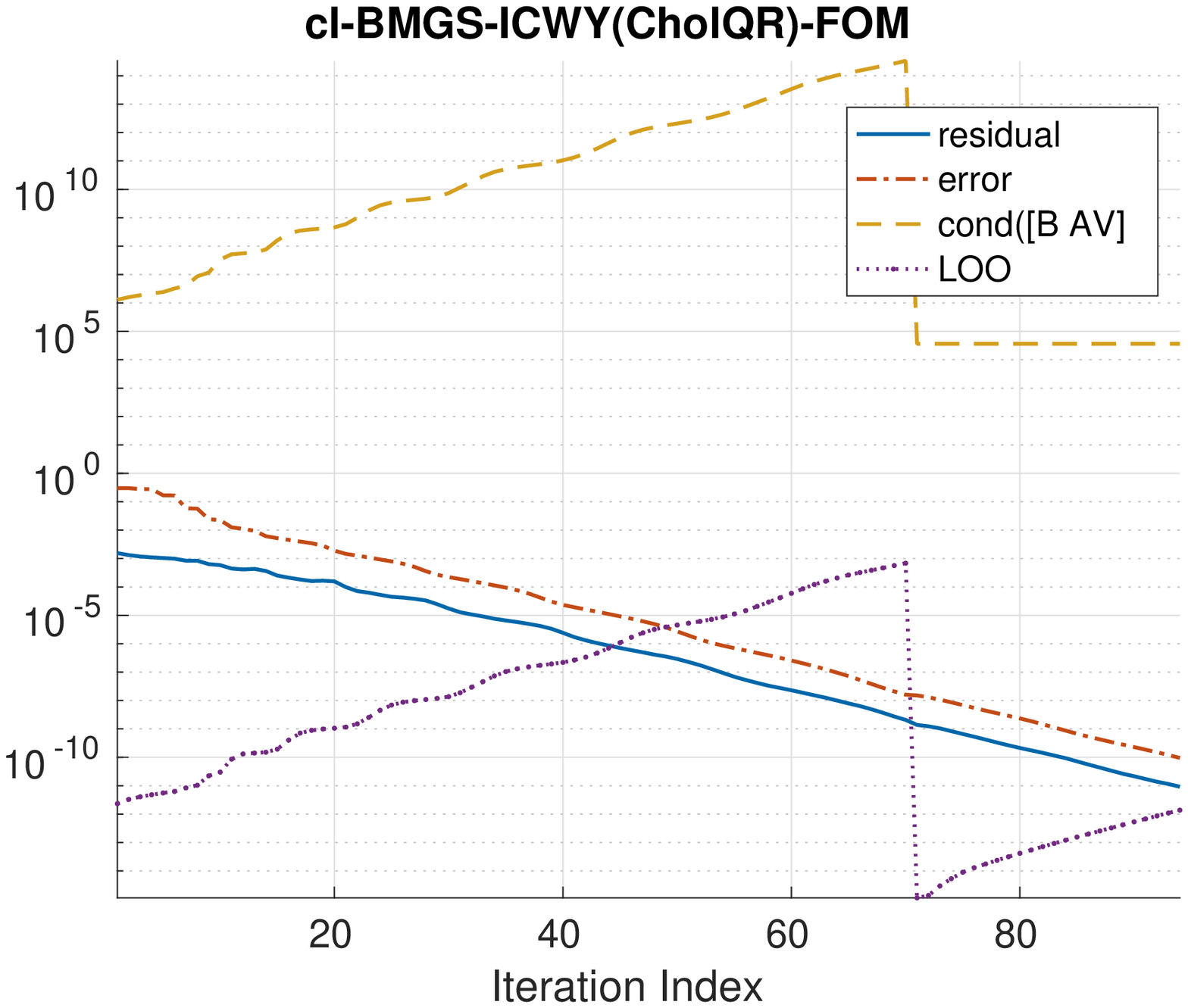}} &
		\resizebox{.44\textwidth}{!}{\includegraphics{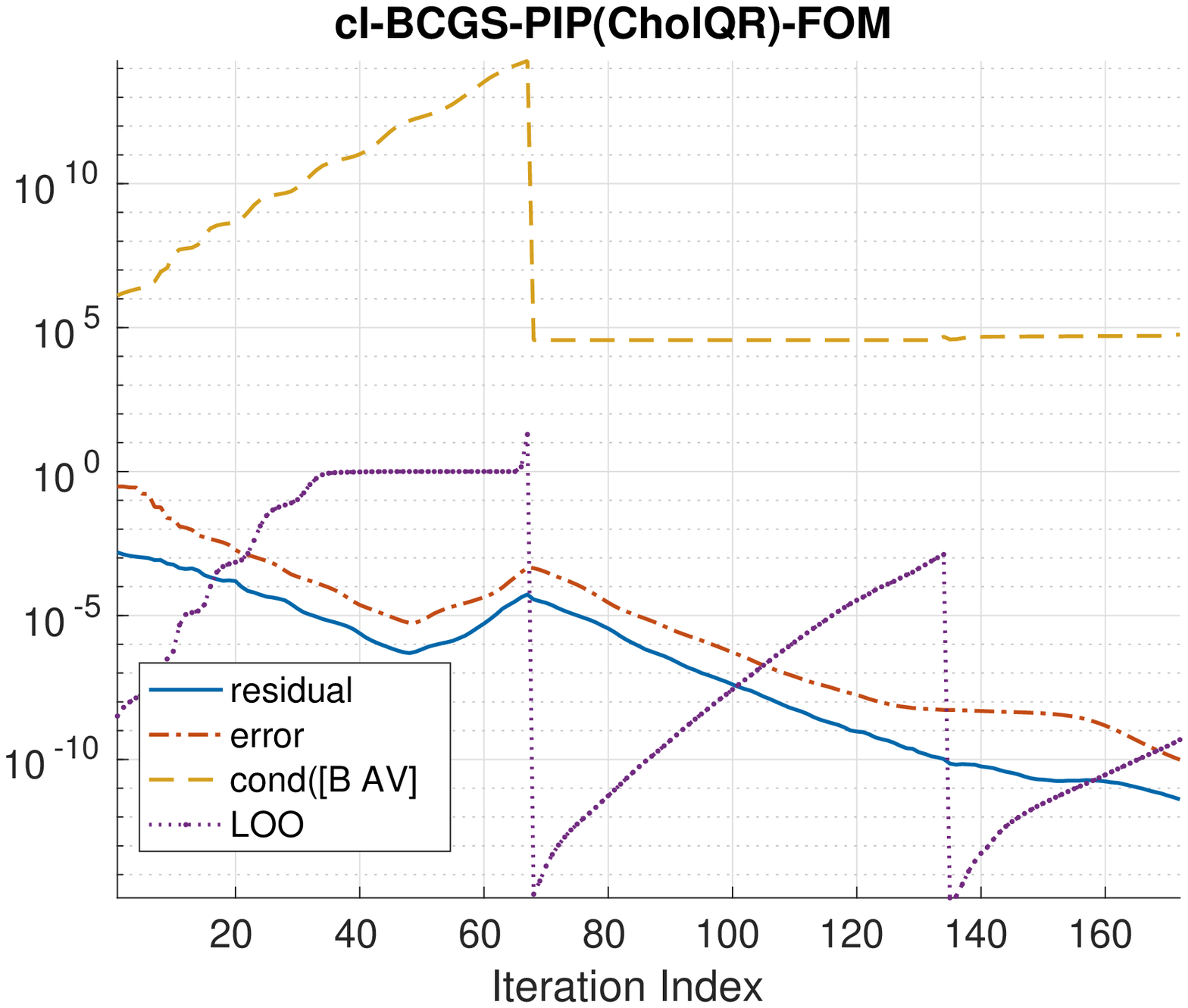}} \\
	\end{tabular}
	\caption{Convergence histories of some classical variants for \tridiag example. \label{fig:tridiag_cl_conv}}
\end{figure}

Figures~\ref{fig:tridiag_gl_conv} and \ref{fig:tridiag_cl_conv} display convergence histories for a subset of the methods in Table~\ref{tab:tridiag}.  The convergence histories for all global \BMGS variants are very similar; we omit \BMGSSVL and \BMGSCWY, as they are visually identical to \BMGSLTS and \BMGSICWY, respectively.  \BMGS is identical to \BMGSSVL and \BMGSLTS and is therefore also omitted.

Both the classical and global variants of \BCGSPIP show the robustness of the adaptive restarting procedure in action.  In the global case, the LOO exceeds $\bigO{10^{-10}}$ and reaches $\bigO{1}$ in \cl-\BCGSPIP.  Despite the loss of orthogonality, restarting allows the methods to recover and eventually converge.  All other low-sync variants remain stable, only restarting once the basis size limit of $m = 70$ has been reached.  Although hardly perceptible, \BMGSICWY does have a slightly worse LOO than that of \BMGSLTS, which can be seen by zooming in on the last few iterations of the global plots in Figure~\ref{fig:tridiag_gl_conv} or of the classical plots in Figure~\ref{fig:tridiag_cl_conv}.

We also note that the residual estimate \eqref{eq:Rm} for all methods follows the same qualitative trend as that of the error.  In the worst case, \cl-\BCGSPIP, the residual is nearly 3 orders of magnitude lower than the error in some places, which could lead to premature convergence.  For all other methods, the difference is between 1 and 2 orders of magnitude.  We would thus recommend setting the residual tolerance a couple orders of magnitude lower in practice, to ensure that the true error is accurate enough.

\subsection{\powersystem} \label{sec:1138_bus}
Now we turn to a slightly more complicated matrix.  The matrix $A$ comes from a power network problem and is real and symmetric positive definite, while \updated{entries of $\vB$ are drawn randomly} from the uniform distribution.  Moreover we apply an incomplete LU (ILU) preconditioner with no fill, using \MATLAB's built-in \texttt{ilu}.

Even with the preconditioner, none of the global methods converges.  We adjusted the thread count to see if it would aid convergence, to no avail.  This is perhaps an extreme case of \cite[Theorem~3.3]{FroLS20}, wherein the global method is much less accurate than the classical method in the first cycle and cannot manage to catch up even after restarting.  A preconditioner better attuned to the structure of the problem may alleviate stagnation for global methods, but we do not explore this here.

\begin{figure}[htbp!]
	\resizebox{.95\textwidth}{!}{\includegraphics{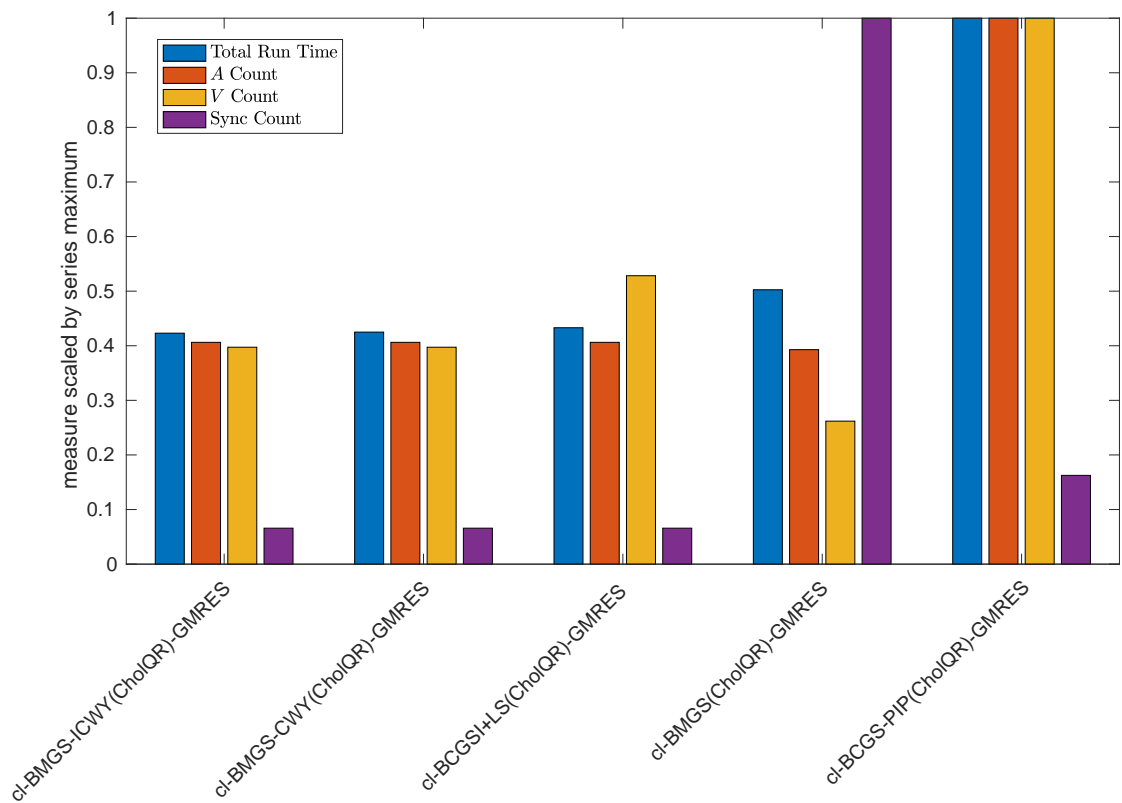}}
	\caption{Performance results for \powersystem example. \label{fig:1138_bus_perf}}
\end{figure}

In Figure~\ref{fig:1138_bus_perf} we see the performance results for the convergent classical methods; more details can be found in Table~\ref{tab:1138_bus}.  Most notably, the one-sync methods \BMGSCWY, \BMGSICWY, and \BCGSIROLS improve over \BMGS only slightly in terms of timings.  \BCGSPIP is much slower, due to a quick loss of orthogonality and need to restart more often.  However, it is clear that sync counts for all one-sync methods are drastically reduced compared to that of \BMGS.

\begin{figure}[htbp!]
	\begin{tabular}{cc}
		\resizebox{.45\textwidth}{!}{\includegraphics{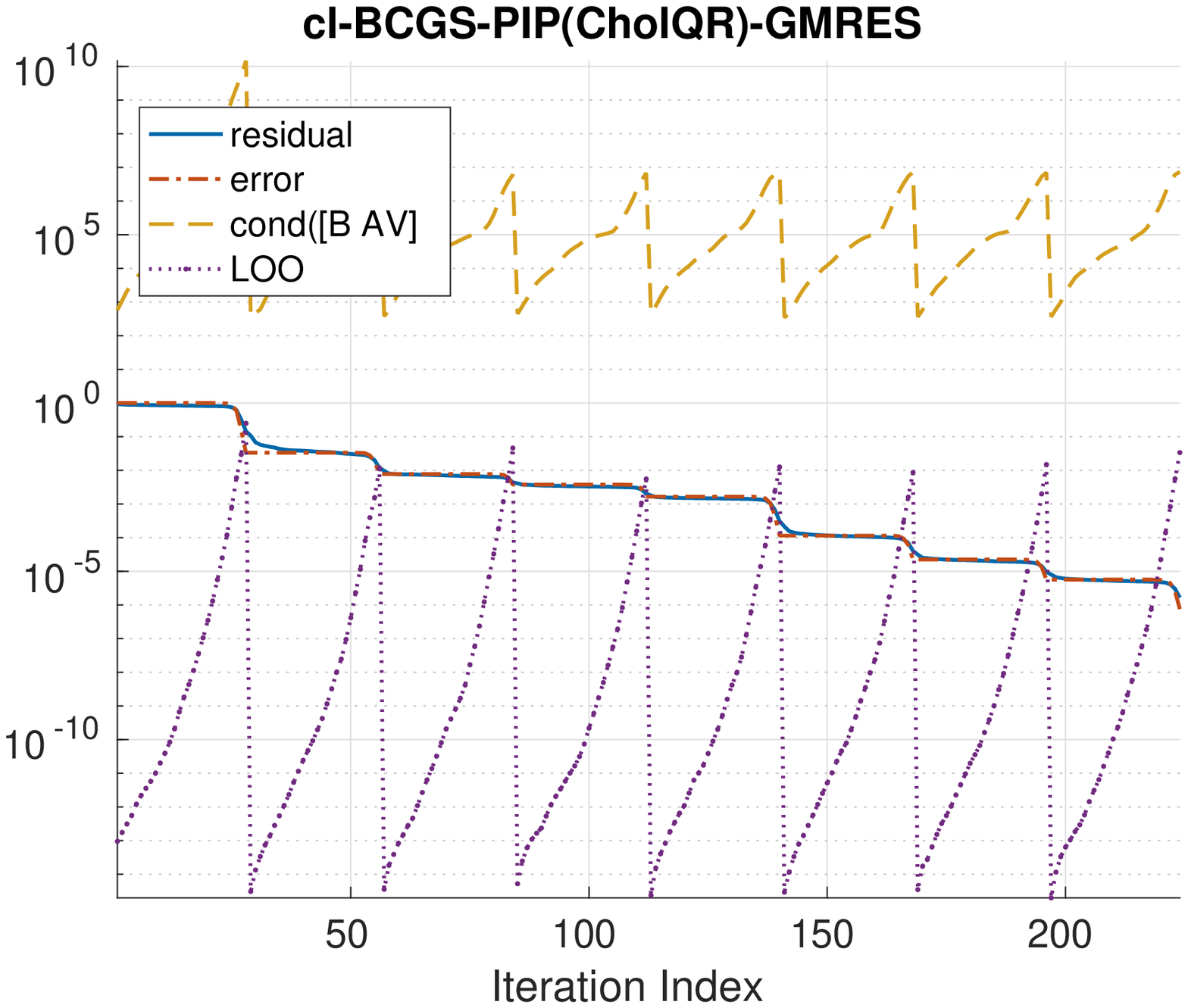}} &
		\resizebox{.45\textwidth}{!}{\includegraphics{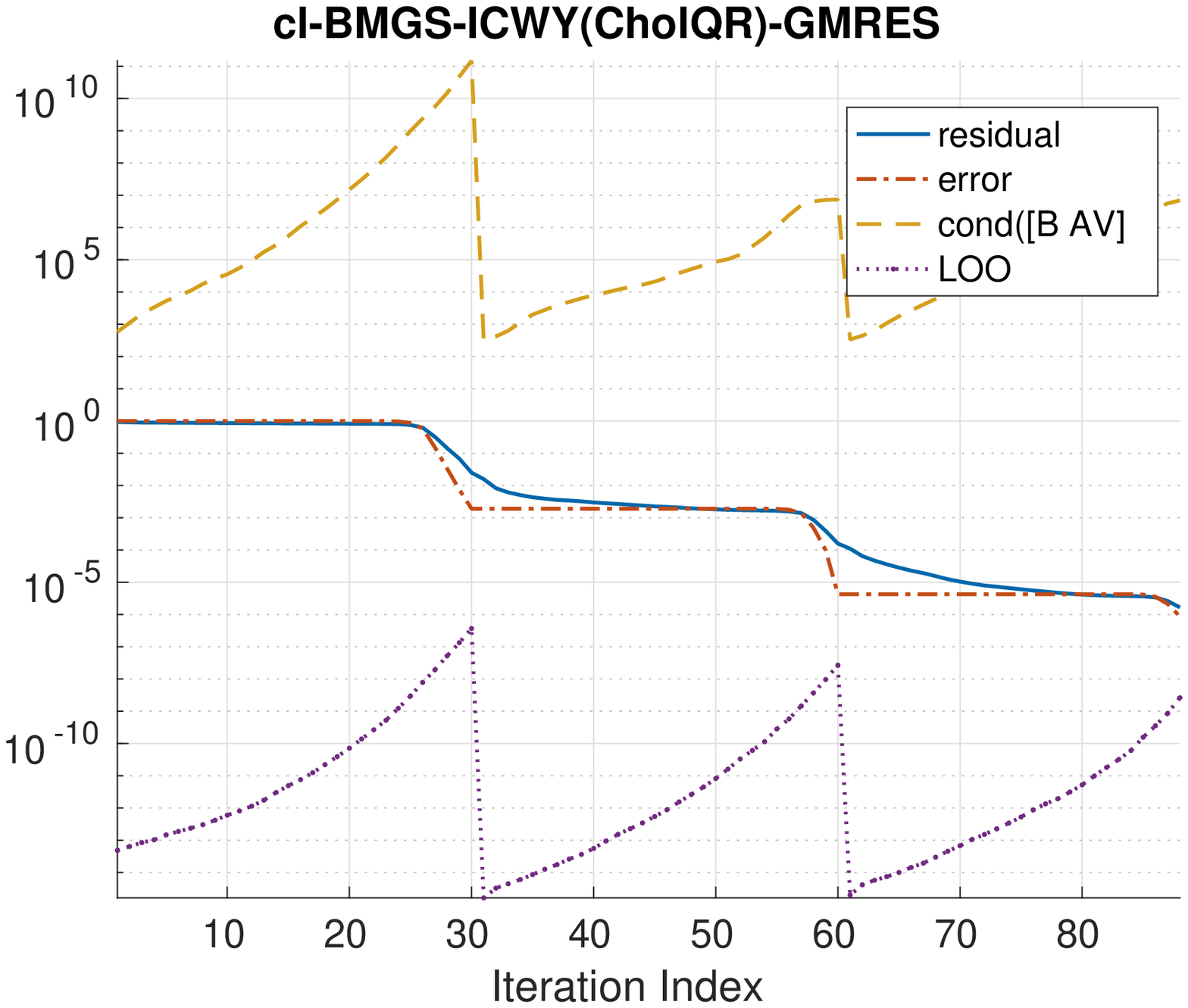}}
	\end{tabular}
	\caption{Subset of convergence histories for \powersystem example. \label{fig:1138_bus_conv}}
\end{figure}

We examine the convergence histories of \cl-\BCGSPIP and \cl-\BMGSICWY more closely in Figure~\ref{fig:1138_bus_conv}.  Although not discernible on the graph, we found that \cl-\BCGSPIP actually restarts every 28 iterations, meaning in the first cycle it encountered a \nan-flag and reduced the maximum basis size to $m = 28$ for all subsequent cycles.  Instability in the first cycle thus hinders \cl-\BCGSPIP greatly.  On the other hand, \BMGSICWY (as well as the other variants) is stable enough to exhaust the entire basis size allowance, which allows for further error reduction in the first cycle.

\subsection{\texttt{circuit\_2}} \label{sec:circuit_2}
The next example comes from a circuit simulation problem.  The matrix $A$ is real but not symmetric or positive definite.  We again apply an ILU preconditioner with no fill.

\begin{figure}[htbp!]
	\resizebox{.95\textwidth}{!}{\includegraphics{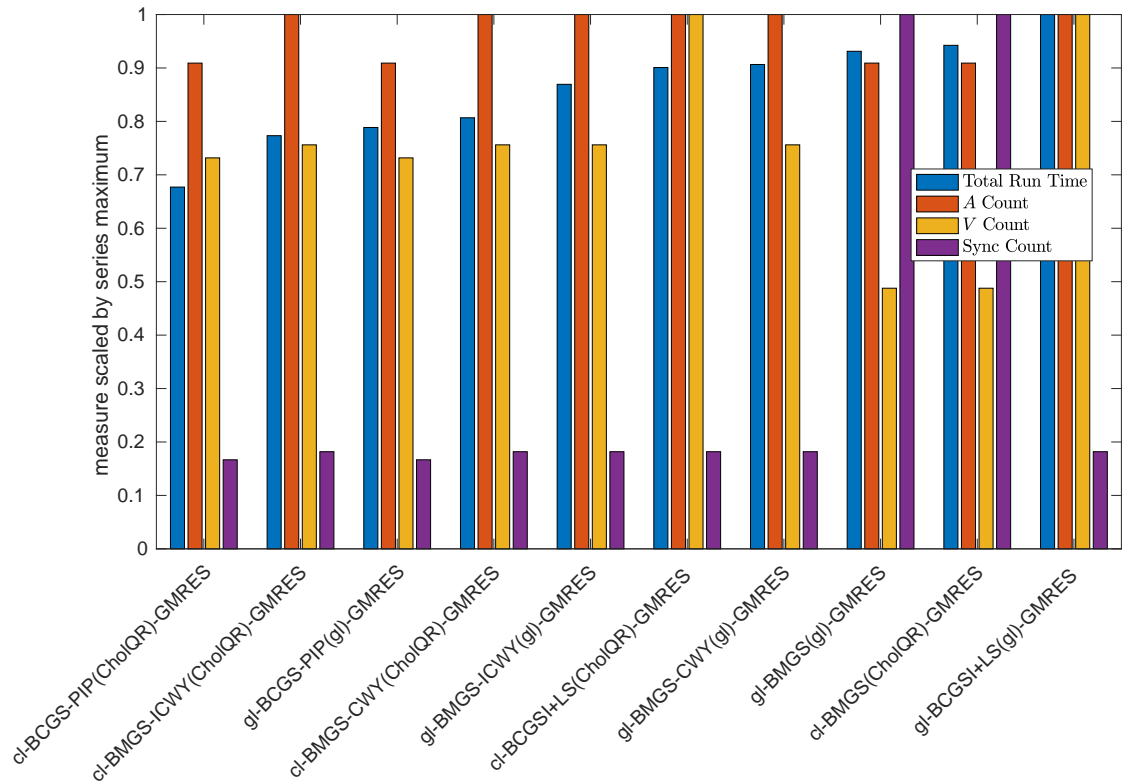}}
	\caption{Performance results for \texttt{circuit\_2} example. \label{fig:circuit_2_perf}}
\end{figure}

All the one-sync classical and global methods converge, and their performance data is presented in Figure~\ref{fig:circuit_2_perf} with further details in Table~\ref{tab:circuit_2}.  In fact, some global methods, like \gl-\BCGSPIP, are even faster than some classical methods, due to the fact that they require the same number of iterations to converge, and therefore fewer floating-point operations.

\begin{figure}[htbp!]
	\begin{tabular}{cc}
		\resizebox{.45\textwidth}{!}{\includegraphics{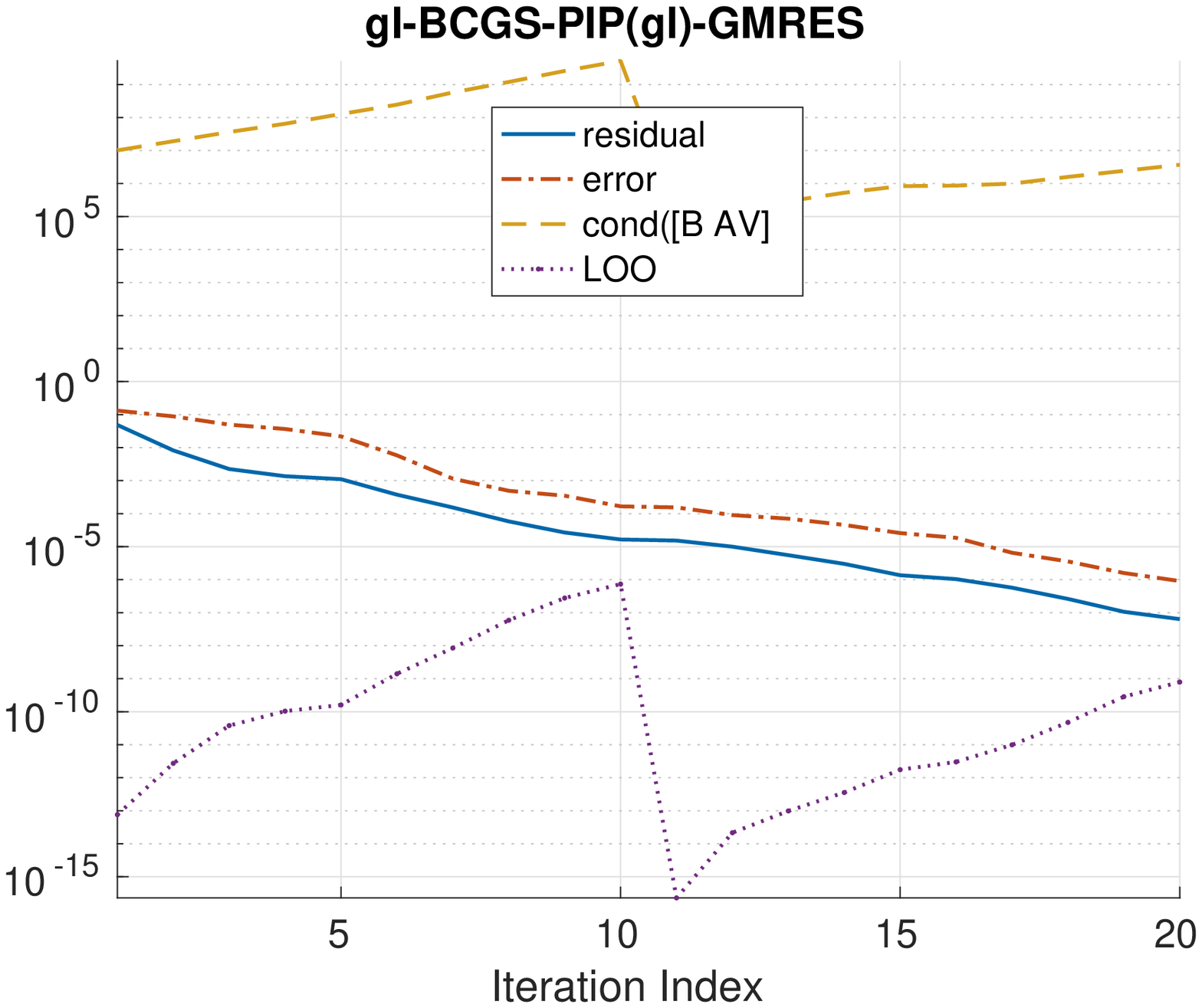}} &
		\resizebox{.45\textwidth}{!}{\includegraphics{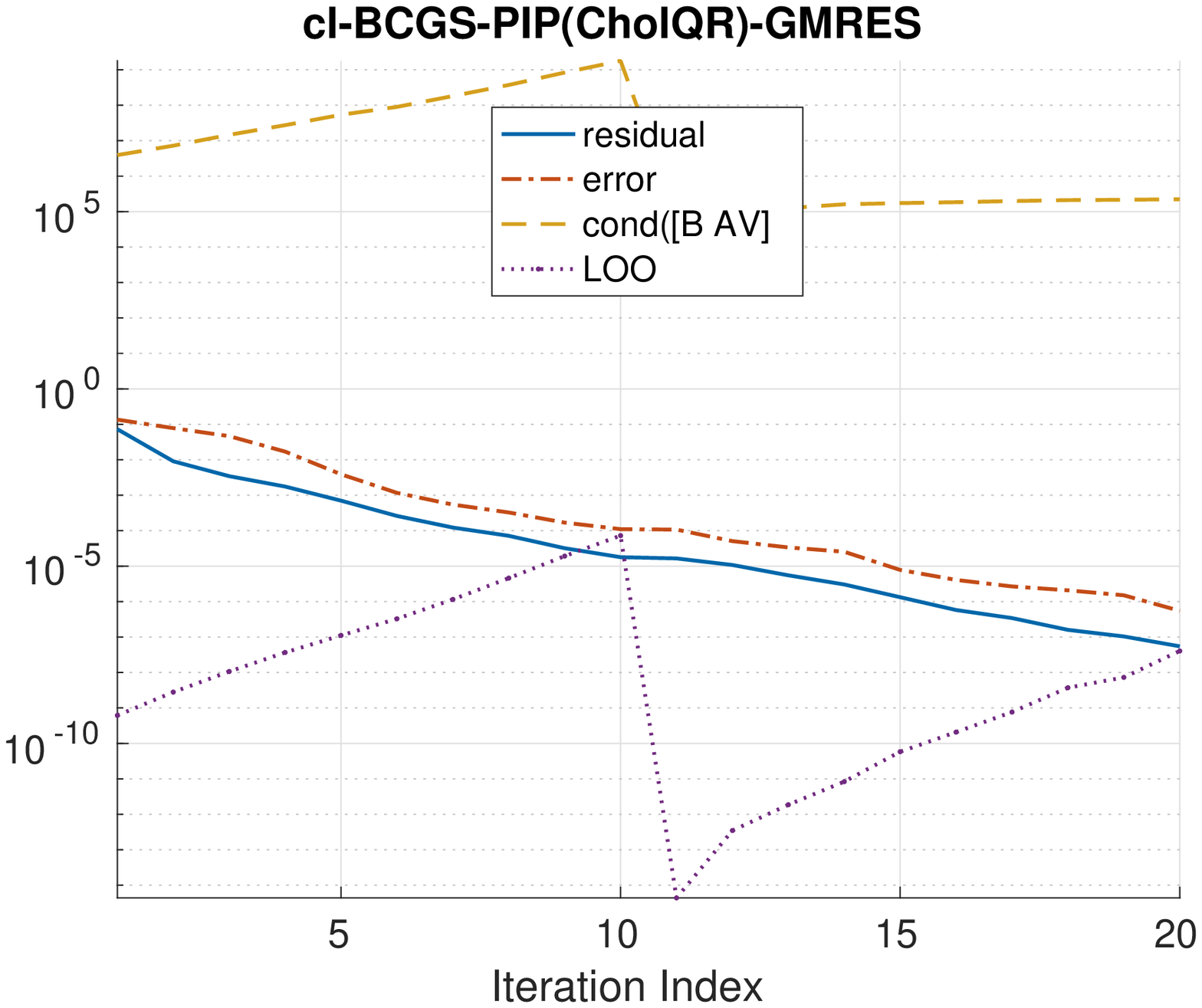}} \\
	\end{tabular}
	\caption{Convergence histories of the \BCGSPIP variants for the \texttt{circuit\_2} example. \label{fig:circuit_2_conv}}
\end{figure}

Figure~\ref{fig:circuit_2_conv} demonstrates how close in accuracy the global and classical \BCGSPIP variants are for this problem.  The global method even has a slightly better LOO, but it should be noted that global LOO is measured according to a different inner product than classical LOO; see Section~\ref{sec:block_inner_prod} and \eqref{eq:loo}.

\subsection{\texttt{rajat03}} \label{sec:rajat03}
Another circuit simulation problem highlights slightly different behavior.  In this case, $A$ is again real but neither symmetric nor positive definite, and we again use an ILU preconditioner with no fill.

\begin{figure}[htbp!]
	\resizebox{.95\textwidth}{!}{\includegraphics{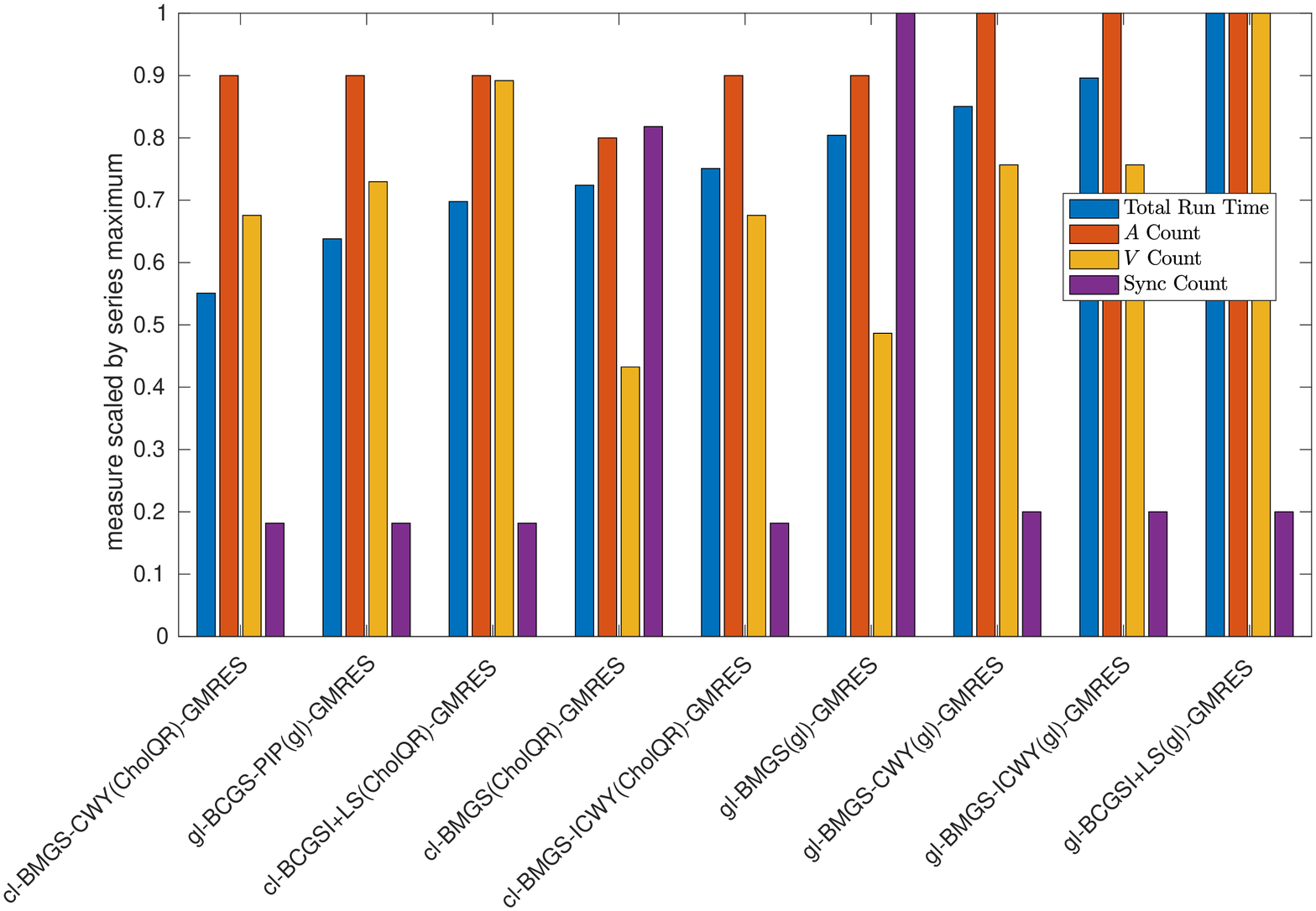}}
	\caption{Performance results for \texttt{rajat03} example. \label{fig:rajat03_perf}}
\end{figure}

Figure~\ref{fig:rajat03_perf} summarizes the performance results, with details given in Table~\ref{tab:rajat03}.  It should be noted right away that \cl-\BCGSPIP fails to converge for this problem, while \gl-\BCGSPIP does not, and takes second place in terms of the timings.  More specifically, \cl-\BCGSPIP encounters a \nan-flag it cannot resolve, which means that every time it reduces the basis size, it cannot avoid a \nan-flag.  However, because global methods do not use Cholesky at all, non-positive definite factors do not pose a problem, unless their trace is numerically zero, which occurs with very low probability.  Otherwise, \cl-\BMGSCWY shows a small improvement over \cl-\BMGS.

\begin{figure}[htbp!]
	\begin{tabular}{cc}
		\resizebox{.45\textwidth}{!}{\includegraphics{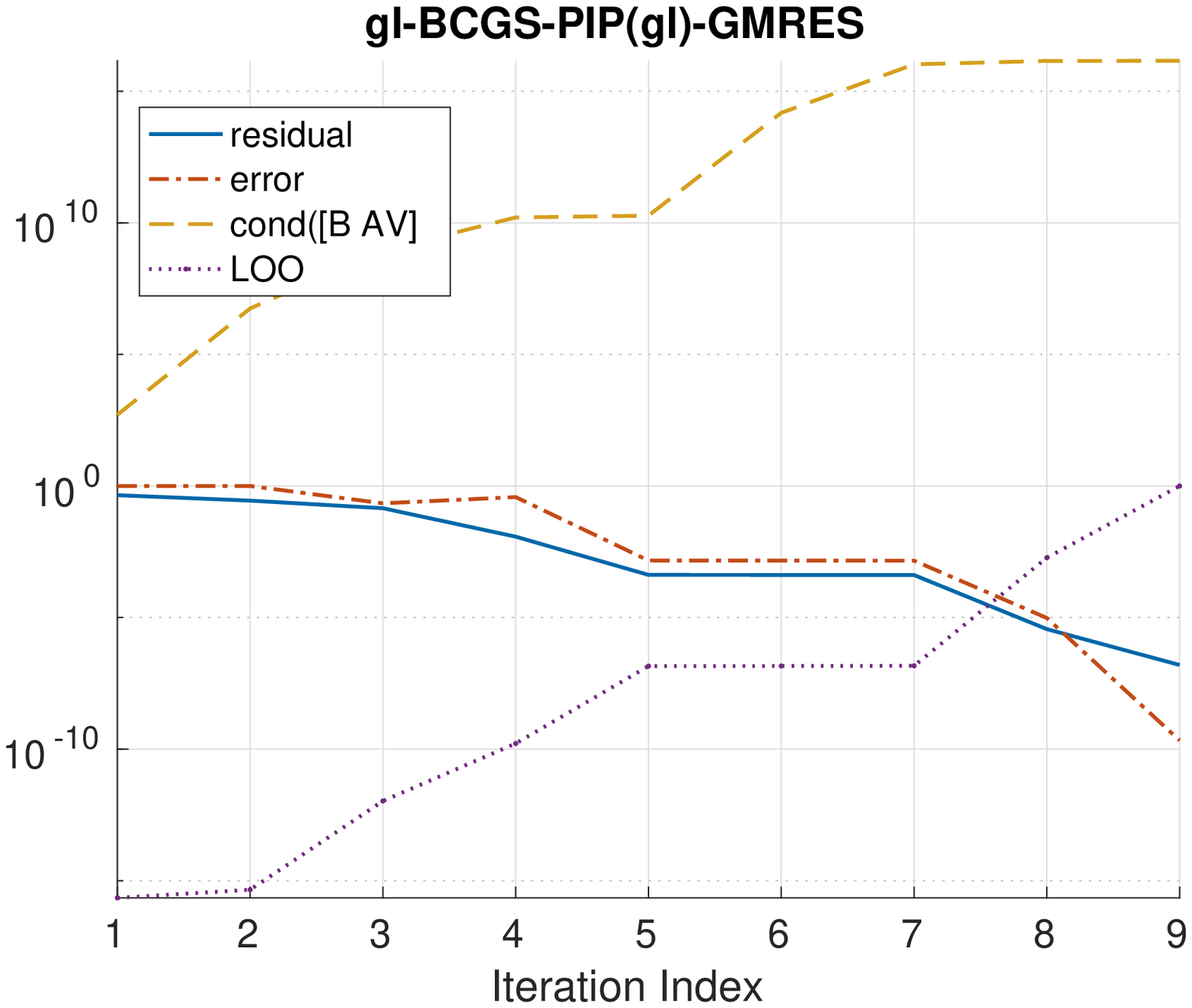}} &
		\resizebox{.45\textwidth}{!}{\includegraphics{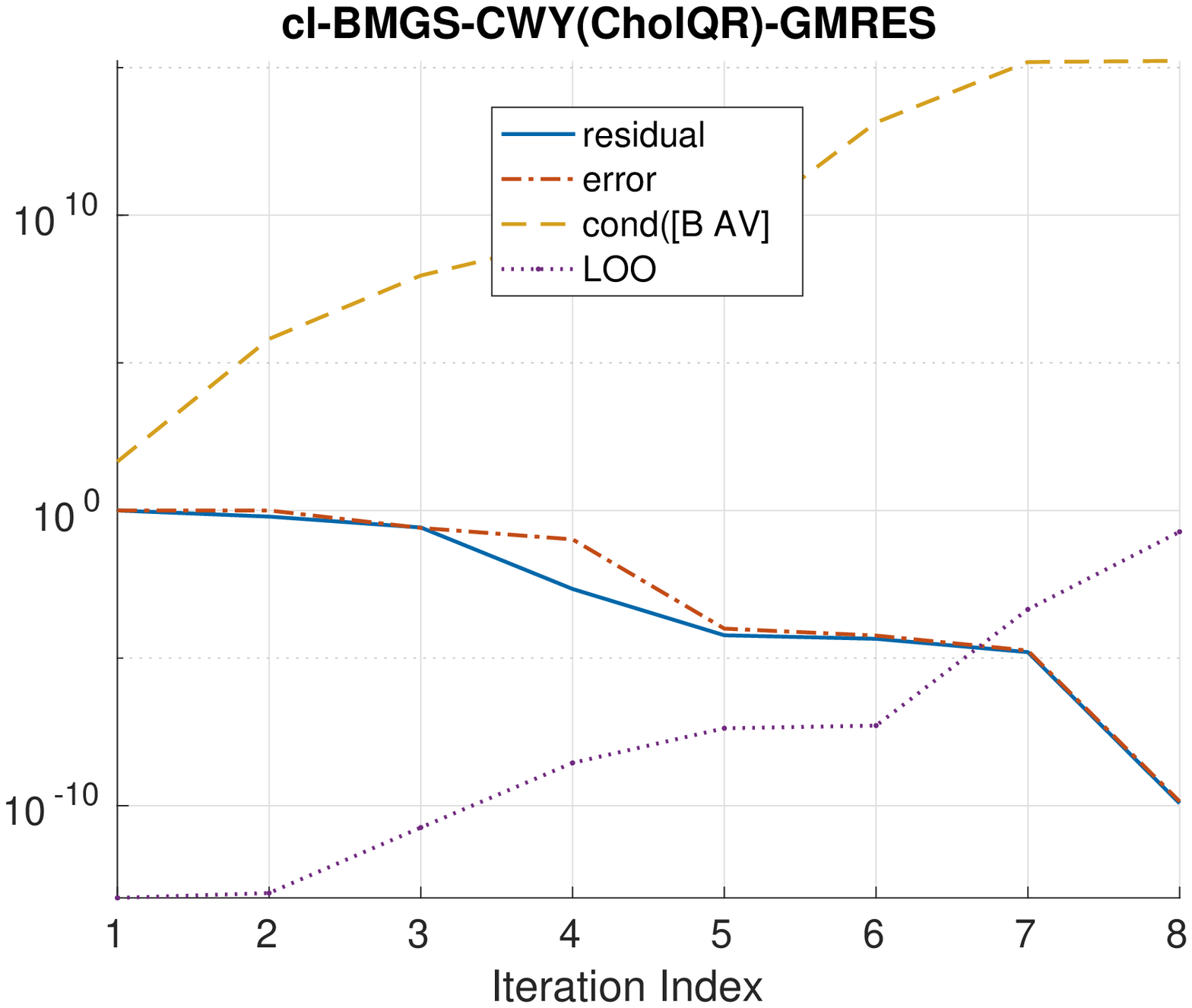}} \\
	\end{tabular}
	\caption{Convergence histories of the two fastest variants for the \texttt{rajat03} example. \label{fig:rajat03_conv}}
\end{figure}

Table~\ref{tab:rajat03} confirms that none of the methods requires restarting despite how high the condition number becomes in later iterations; see also Figure~\ref{fig:rajat03_conv}.  It is again interesting to see how close the error and residual plots are between the global and classical methods.  In fact, the residual for the global method underestimates convergence by a couple orders of magnitude.

\subsection{\texttt{Kaufhold}} \label{sec:kaufhold}
This example treats a nearly numerically singular matrix with an extremely high condition number.  Also notable, the norm of $A$ is nearly $\bigO{10^{15}}$.  The matrix is real, but neither symmetric nor positive definite, and it was designed to trigger a bug in Gaussian elimination in a 2002 version of \MATLAB.  We again apply an ILU preconditioner with no fill.

\begin{figure}[htbp!]
	\resizebox{.95\textwidth}{!}{\includegraphics{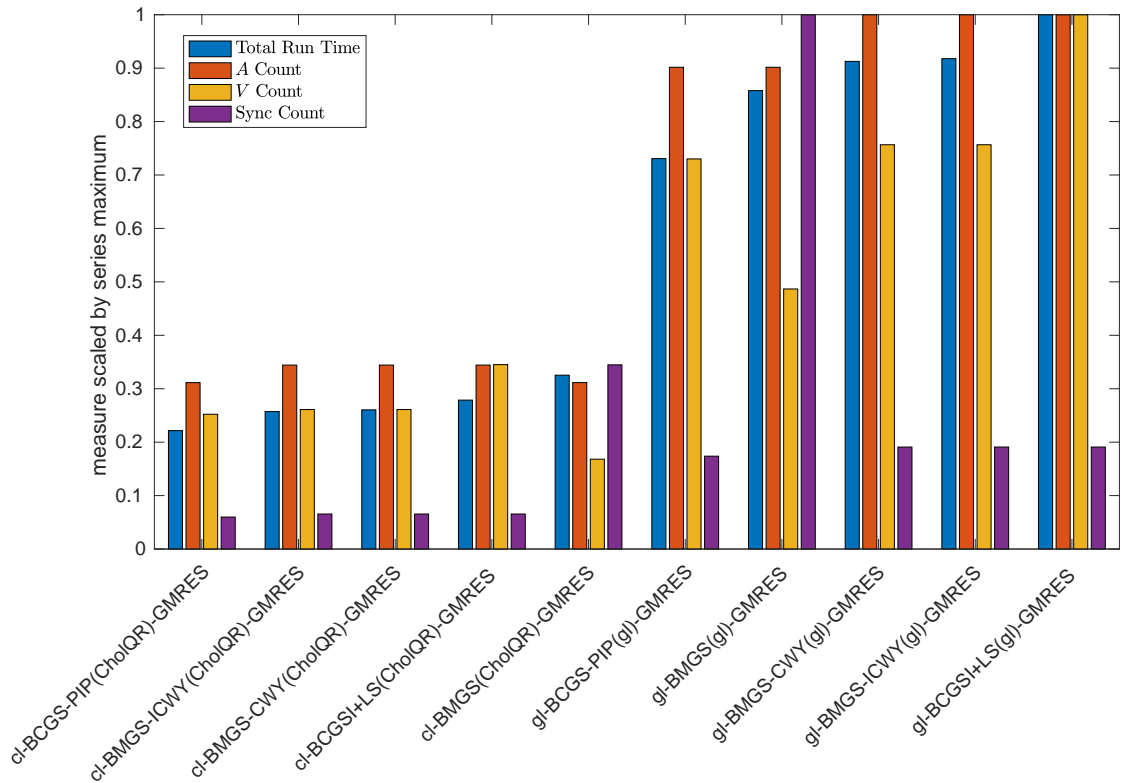}}
	\caption{Performance results for \texttt{Kaufhold} example. \label{fig:Kaufhold_perf}}
\end{figure}

\begin{figure}[htbp!]
	\begin{tabular}{cc}
		\resizebox{.45\textwidth}{!}{\includegraphics{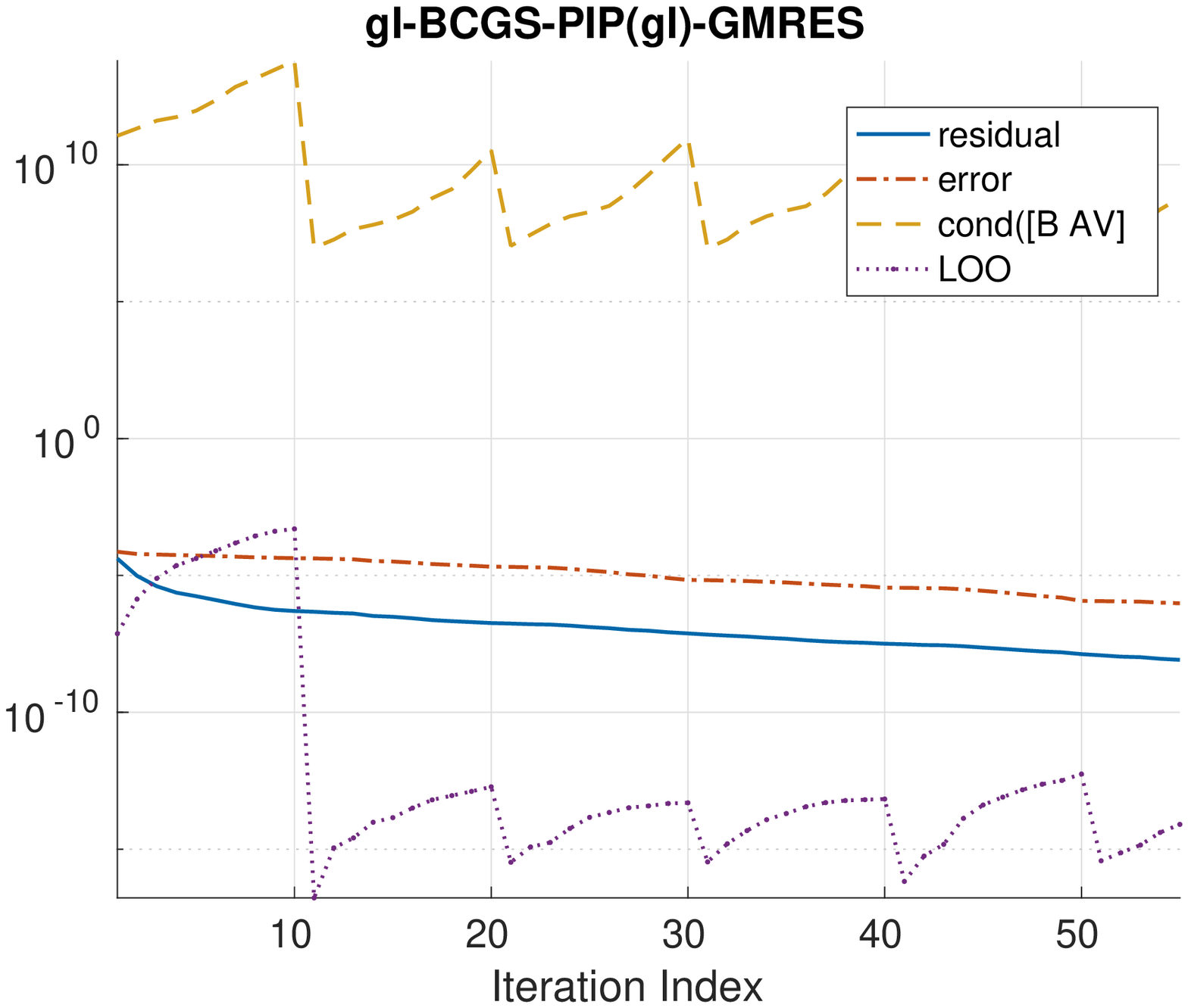}} &
		\resizebox{.45\textwidth}{!}{\includegraphics{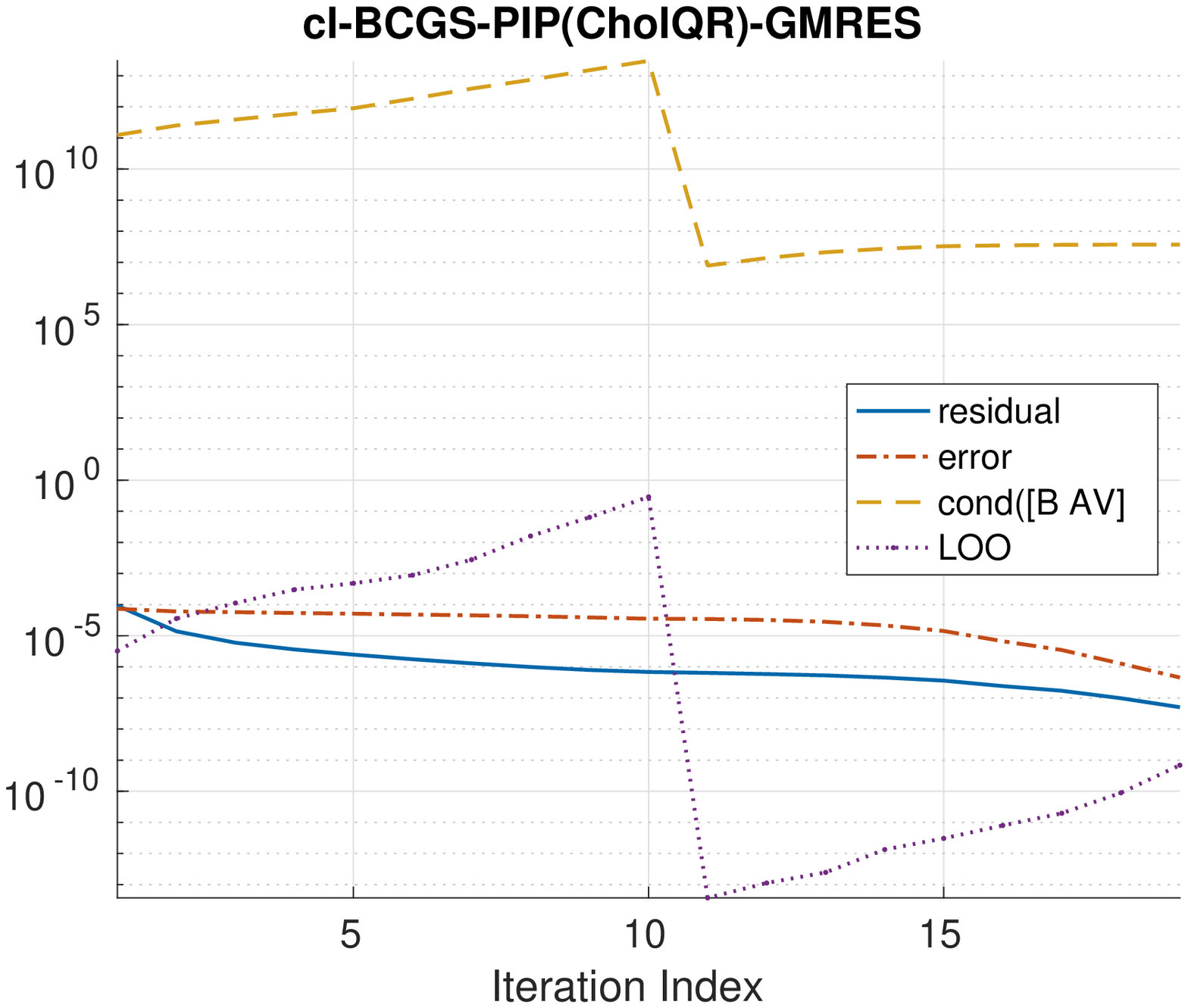}} \\
	\end{tabular}
	\caption{Convergence histories of the \BCGSPIP variants for the \texttt{Kaufhold} example. \label{fig:Kaufhold_conv}}
\end{figure}

Figure~\ref{fig:Kaufhold_perf} shows \cl-\BCGSPIP to be the fastest of the classical one-sync methods, but the improvement over \cl-\BMGS is small.  The global methods are all much slower.  A look at the convergence histories in Figure~\ref{fig:Kaufhold_conv} shows a stubborn error curve despite significant progress in the initial iterations.  For both \BCGSPIP methods the LOO is moderately high in the first cycle, matching the high condition numbers, but the situation is not bad enough to trigger a \nan-flag, and the LOO drops after restarting.

\subsection{\texttt{t2d\_q9}} \label{sec:t2d_q9}
We now examine a nonlinear diffusion problem, specifically a biquadratic mesh of a temperature field.  The matrix $A$ is real but not symmetric or positive definite, and we again use an ILU preconditioner with no fill.

\begin{figure}[htbp!]
	\resizebox{.95\textwidth}{!}{\includegraphics{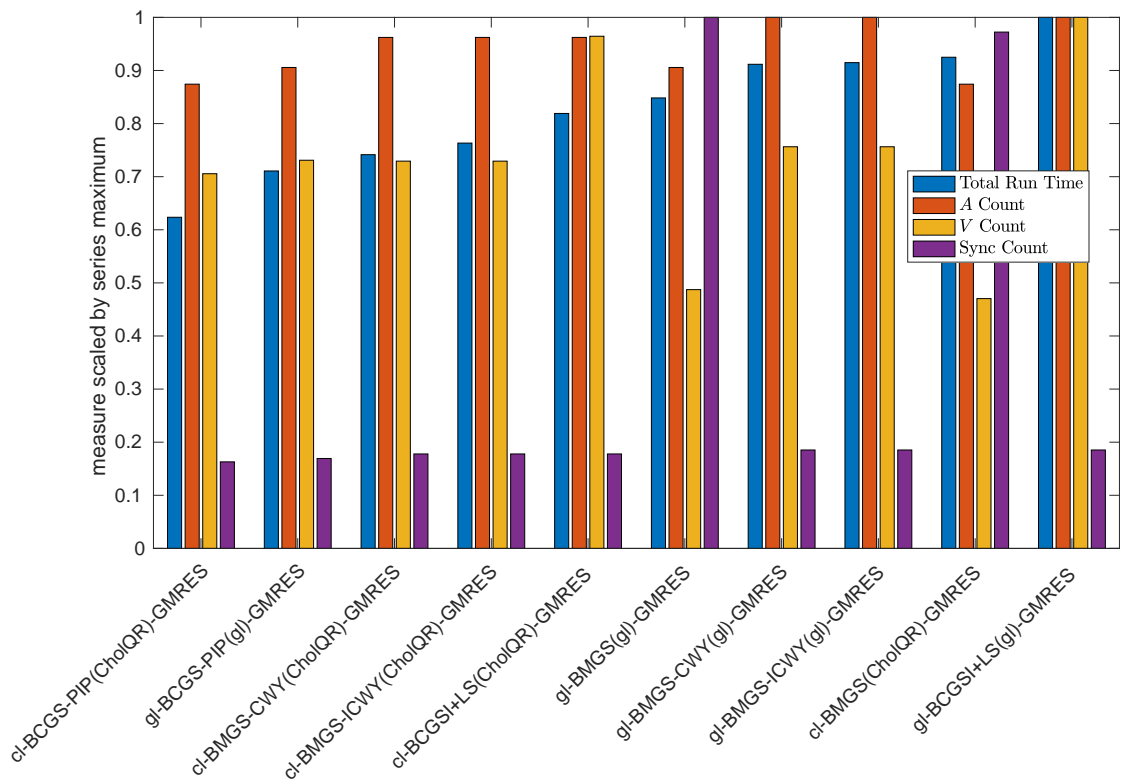}}
	\caption{Performance results for \texttt{t2d\_q9} example. \label{fig:t2d_q9_perf}}
\end{figure}

Figure~\ref{fig:t2d_q9_perf} shows that both \BCGSPIP are the fastest overall, with \cl-\BMGS in second-to-last place; see Table~\ref{tab:t2d_q9} for more details.  Interestingly, even \gl-\BMGS is faster than \cl-\BMGS in this scenario.

Both \BCGSIROLS variants are rather slow in this example.  Despite having just one sync per iteration, \BCGSIROLS does generally have a higher complexity than its one-sync counterparts, which manifests here as a disadvantage.

\begin{figure}[htbp!]
	\begin{tabular}{cc}
		\resizebox{.45\textwidth}{!}{\includegraphics{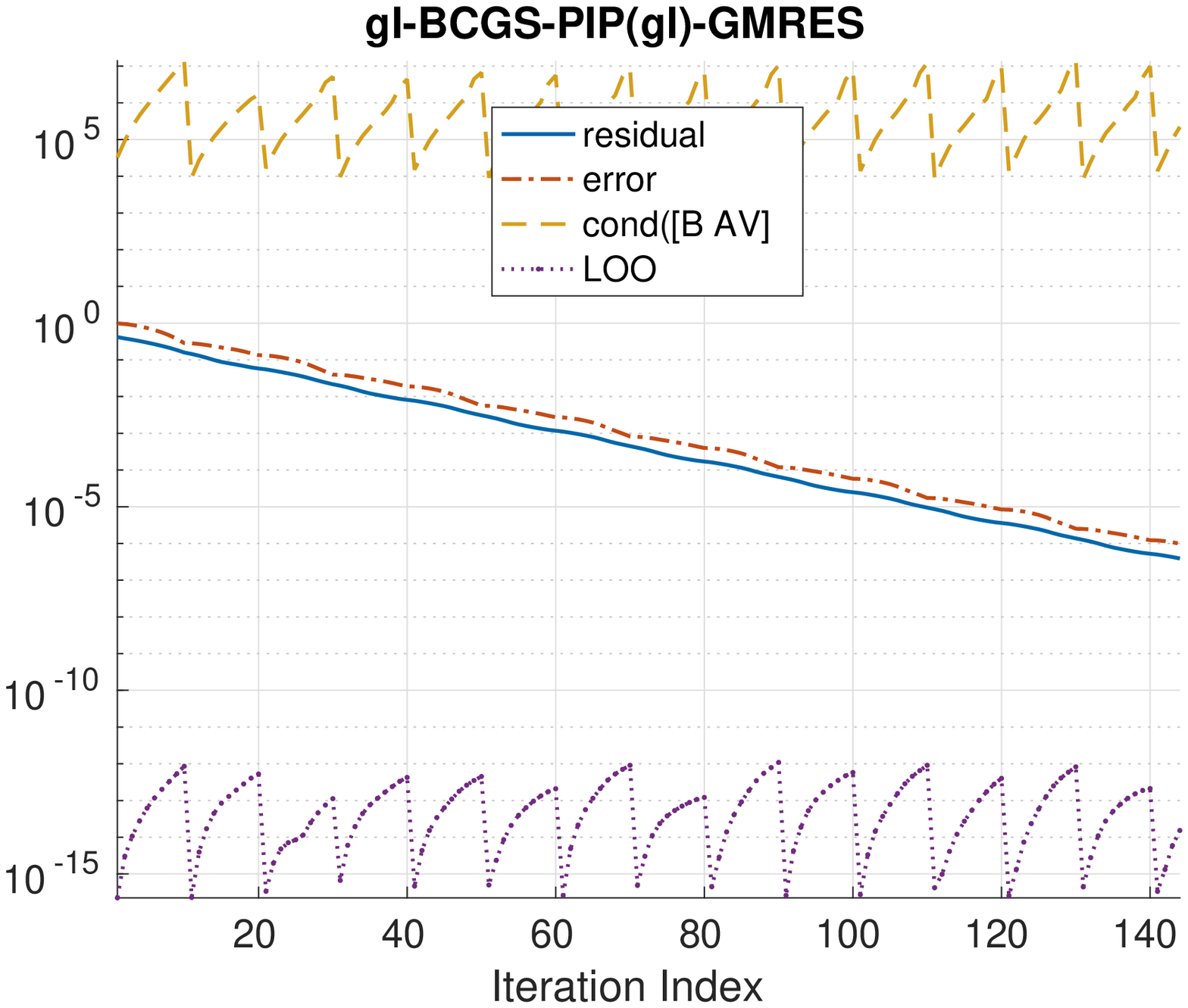}} &
		\resizebox{.45\textwidth}{!}{\includegraphics{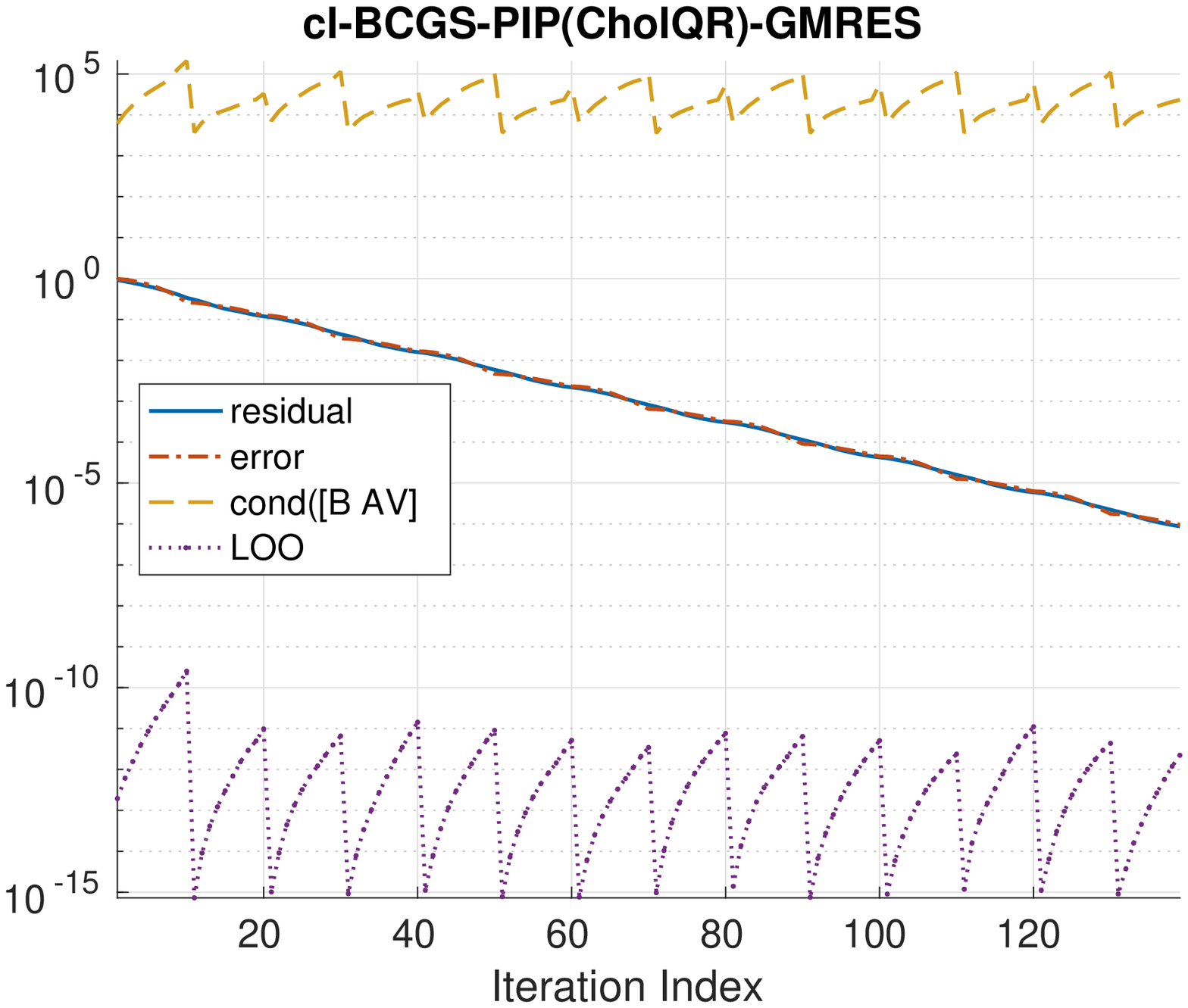}} \\
	\end{tabular}
	\caption{Convergence histories of the \BCGSPIP variants for the \texttt{t2d\_q9} example. \label{fig:t2d_q9_conv}}
\end{figure}

The convergence behavior for the \BCGSPIP variants is given in Figure~\ref{fig:t2d_q9_conv}.  Here we see that despite the global condition number having a high variation relative to the classical method, the global LOO is overall much less.  This phenomenon is not unique to this example, however, it just happens to be more noticeable.

\subsection{\lapl} \label{sec:lapl_2d}
Our last problem is taken directly from \cite[Section~5.4]{FroLS17}, a discretized two-dimensional Laplacian matrix.  $A$ is thus banded, real, and symmetric positive definite.  We do not apply a preconditioner, and look at all skeletons considered in the text.

\begin{figure}[htbp!]
	\resizebox{.95\textwidth}{!}{\includegraphics{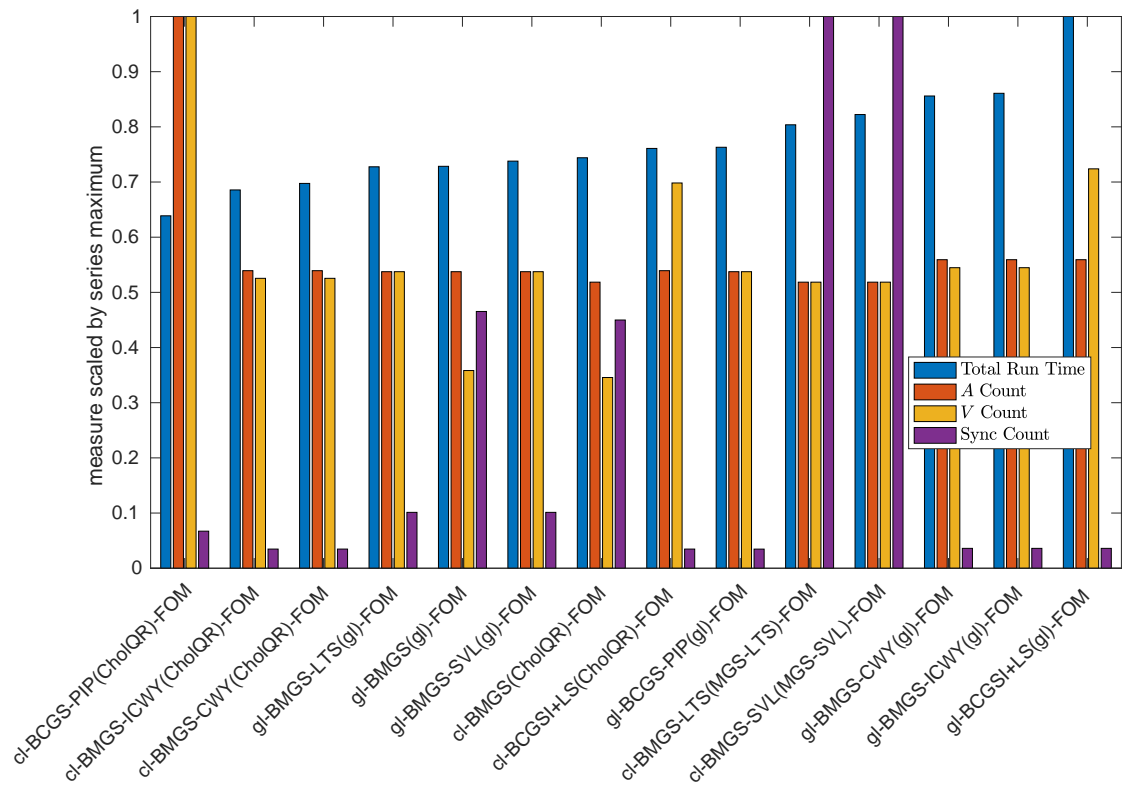}}
	\caption{Performance results for \lapl example. \label{fig:lapl_2d_perf}}
\end{figure}

Figure~\ref{fig:lapl_2d_perf} shows the performance results; more details can be found in Table~\ref{tab:lapl_2d}.  All one-sync classical methods except for \cl-\BCGSIROLS beat \cl-\BMGS, along with a number of global methods.  The slowest classical methods are the three-sync ones, and some one-sync global methods follow behind.  The fastest method, \cl-\BCGSPIP also happens to have the highest $A$ count and applications of $\bVV_k$, due to its high number of restarts.  Both \cl-\BMGSCWY and \cl-\BMGSICWY, however, have fewer sync counts, as well as $A$ counts and $\bVV_k$ counts, and are very close in terms of timings.

The methods with the highest sync counts are \cl-\BMGSSVL and \cl-\BMGSLTS.  The reason is that they cannot use \CholQR as a muscle,\footnote{Strictly speaking, they can use whatever muscle they are programmed to use, but \BMGSSVL requires \MGSSVL to be stable; see \cite{Bar19, CarLRetal22}.} and this problem requires many iterations to converge.  \LSBA is written to count sync points within the muscles as well, and with \MGSSVL and \MGSLTS each contributing $1 + 3s$ per call, the total number of sync points eventually passes that of \cl-\BMGS, which can use a communication-light muscle like \CholQR.

\section{Conclusions and outlook} \label{sec:conclusion}
Stability bounds and floating-point analysis are challenging to work out rigorously, and it is therefore simultaneously important to search for counterexamples and edge cases while trying to prove conjectured bounds.  In general, rigorous loss of orthogonality and backward error bounds for all these methods could lead to new insights and improvements in the quest for a reliable, scalable Krylov subspace solver.  Our flexible benchmarking tool can aid in that process, and it can easily be extended to accommodate new algorithm \updated{configurations, test cases, and measures}.

At the same time, low-sync block Arnoldi algorithms with adaptive restarting are clearly already useful and robust enough for a wide variety of problems, especially where $A$ is reasonably conditioned and memory limitations cap basis sizes.  In every benchmark, we have observed that at least one low-sync method outperformed both the classical and global \BMGS-based Arnoldi methods.  \updated{More research is needed to determine which low-sync skeletons are best for which problems and architectures, particularly computational models that account not only for operation counts but also for performance variations relative to block size \cite{Bir15, BomHS22, ParSS16}.  Most likely the best configuration allows for switching between skeletons and muscles depending on convergence behavior.}

For scenarios where the basic adaptive restarting procedure is not sufficient to rescue convergence, it might be possible to improve the heuristics with a cheap estimate of the loss of orthogonality computed, e.g., a randomized sketched inner product \cite{BalGr22}.  With such a cheap estimate, we could not only decrease the basis size when there are problems, but increase it again in later cycles.  Randomized algorithms themselves are known to reduce communication, and a thorough comparison and combination of the methods proposed here and in \cite{BalGr22} could lead to powerful Krylov subspace method well suited for exascale architectures.

\updated{Global methods are unfortunately less promising.  They are almost always slower than even the slowest classical method, due to requiring more cycles, and thus operator calls and sync points, to converge.  However, the benchmarks do suggest that, in cases with a good preconditioner known to guarantee convergence in a few iterations, global methods may become competitive again, especially in single-node or ``laptop" applications, where their reduced computational intensity per iteration is favorable.}

\section{Declarations}

\paragraph{Ethical Approval and Consent to participate}
The author certifies that this manuscript has been submitted to only one journal at this time, that the work is original, and that the results are not fabricated or skewed.  The work is entirely the author's own, and to the best of the author's ability, the work is complete in its own right and without error or misappropriation.

\paragraph{Consent for publication}
As the sole author, K.~Lund provides consent for publication.

\paragraph{Human and Animal Ethics}
Not applicable.

\paragraph{Availability of supporting data}
All code and scripts to reproduce plots can be found at \url{https://gitlab.mpi-magdeburg.mpg.de/lund/low-sync-block-arnoldi}.

\paragraph{Competing interests}
The author has no relevant financial or non-financial interests to disclose.

\paragraph{Funding}
K.~Lund is a contracted employee of Max Planck Institute for Dynamics of Complex Technical Systems and did not receive any additional funding to support this project.

\paragraph{Authors' contributions}
K.~Lund is the sole author of the manuscript and associated code.

\paragraph{Acknowledgments} The author is indebted to St\'{e}phane Gaudreault, Teodor Nikolov, and Erin Carson for stimulating discussions that inspired this work.  The author is also grateful to Jens Saak and Martin K\"ohler for answering questions about the Mechthild cluster and multithreading in MATLAB \updated{and to two anonymous reviewers for their constructive feedback.}


\addcontentsline{toc}{section}{References}
\bibliography{low_sync_block_arnoldi}
\bibliographystyle{plainurl}


\begin{appendix}
	\section{Raw data from tests}
	\updated{
		A subset of raw data corresponding to the performance plots in Section~\ref{sec:examples} is provided below.  Many headers are abbreviated for space reasons: ``Accel." refers to ``acceleration" or ``speed-up"; ``Ct." refers to ``Count"; and ``Iter." refers to ``Iteration."
	}
	
	\begin{table}[htbp!]
		\begin{center}
			\begin{tabular}{c|c|c|c|c|c|c|c}
				\toprule
				\thead{Configuration}	& \thead{Time\\(s)}	& \thead{\%\\Accel.}	& \thead{Cycle\\Ct.}	& \thead{Iter.\\Ct.}	& \thead{$A$\\Ct.} & \thead{$V$\\Ct.}	& \thead{Sync\\Ct.} \\
				\midrule
				\gl-$\BMGS\circ \gl$-\fom      & 2.20e+00       & 0.00       & 9          & 621          & 621          & 1242          & 22401          \\\hline
				\gl-$\BMGSSVL\circ \gl$-\fom      & 2.05e+00       & 6.59       & 9          & 621          & 621          & 1863          & 1872          \\\hline
				\gl-$\BMGSLTS\circ \gl$-\fom      & 2.03e+00       & 7.90       & 9          & 621          & 621          & 1863          & 1872          \\\hline
				\gl-$\BCGSIROLS\circ \gl$-\fom      & 1.84e+00       & 16.53       & 9          & 621          & 630          & 2493          & 639          \\\hline
				\gl-$\BMGSCWY\circ \gl$-\fom      & 1.58e+00       & 28.26       & 9          & 621          & 630          & 1872          & 639          \\\hline
				\gl-$\BMGSICWY\circ \gl$-\fom      & 1.48e+00       & 32.86       & 9          & 621          & 630          & 1872          & 639          \\\hline
				\gl-$\BCGSPIP\circ \gl$-\fom      & 1.37e+00       & 37.79       & 9          & 621          & 621          & 1863          & 630          \\\hline
				\cl-$\BCGSPIP\circ \CholQR$-\fom      & 5.33e-01       & 75.77       & 3          & 172          & 172          & 516          & 175          \\\hline
				\cl-$\BMGS\circ \CholQR$-\fom      & 4.08e-01       & 81.44       & 2          & 94          & 94          & 188          & 2881          \\\hline
				\cl-$\BMGSCWY\circ \CholQR$-\fom      & 2.80e-01       & 87.28       & 2          & 94          & 96          & 284          & 98          \\\hline
				\cl-$\BMGSSVL\circ \MGSSVL$-\fom      & 2.79e-01       & 87.32       & 2          & 96          & 96          & 288          & 584          \\\hline
				\cl-$\BMGSLTS\circ \MGSLTS$-\fom      & 2.69e-01       & 87.77       & 2          & 96          & 96          & 288          & 584          \\\hline
				\cl-$\BCGSIROLS\circ \CholQR$-\fom      & 2.56e-01       & 88.35       & 2          & 94          & 96          & 378          & 98          \\\hline
				\cl-$\BMGSICWY\circ \CholQR$-\fom      & 2.23e-01       & 89.87       & 2          & 94          & 96          & 284          & 98   
			\end{tabular}
			\caption{Results from \tridiag example. \label{tab:tridiag}}
		\end{center}
	\end{table}
	
	\begin{table}[htbp!]
		\begin{center}
			\begin{tabular}{c|c|c|c|c|c|c|c}
				\toprule
				\thead{Configuration}	& \thead{Time\\(s)}	& \thead{\%\\Accel.}	& \thead{Cycle\\Ct.}	& \thead{Iter.\\Ct.}	& \thead{$A$\\Ct.} & \thead{$V$\\Ct.}	& \thead{Sync\\Ct.} \\
				\midrule
				\cl-$\BCGSPIP\circ\CholQR$-\gmres      & 1.84e+00       & 0.00       & 8          & 224          & 224          & 672          & 232          \\\hline
				\cl-$\BMGS\circ\CholQR$-\gmres      & 9.25e-01       & 49.75       & 3          & 88          & 88          & 176          & 1427          \\\hline
				\cl-$\BCGSIROLS\circ\CholQR$-\gmres      & 7.97e-01       & 56.70       & 3          & 88          & 91          & 355          & 94          \\\hline
				\cl-$\BMGSCWY\circ\CholQR$-\gmres      & 7.82e-01       & 57.50       & 3          & 88          & 91          & 267          & 94          \\\hline
				\cl-$\BMGSICWY\circ\CholQR$-\gmres      & 7.78e-01       & 57.70       & 3          & 88          & 91          & 267          & 94          
			\end{tabular}
		\end{center}
		\caption{Results from \powersystem example. \label{tab:1138_bus}}
	\end{table}
	
	\begin{table}[htbp!]
		\begin{center}
			\begin{tabular}{c|c|c|c|c|c|c|c}
				\toprule
				\thead{Configuration}	& \thead{Time\\(s)}	& \thead{\%\\Accel.}	& \thead{Cycle\\Ct.}	& \thead{Iter.\\Ct.}	& \thead{$A$\\Ct.} & \thead{$V$\\Ct.}	& \thead{Sync\\Ct.} \\
				\midrule
				\gl-$\BCGSIROLS\circ\gl$-\gmres      & 1.26e-01       & 0.00       & 2          & 20          & 22          & 82          & 24          \\\hline
				\cl-$\BMGS\circ\CholQR$-\gmres      & 1.19e-01       & 5.76       & 2          & 20          & 20          & 40          & 132          \\\hline
				\gl-$\BMGS\circ\gl$-\gmres      & 1.18e-01       & 6.88       & 2          & 20          & 20          & 40          & 132          \\\hline
				\gl-$\BMGSCWY\circ\gl$-\gmres      & 1.15e-01       & 9.35       & 2          & 20          & 22          & 62          & 24          \\\hline
				\cl-$\BCGSIROLS\circ\CholQR$-\gmres      & 1.14e-01       & 9.92       & 2          & 20          & 22          & 82          & 24          \\\hline
				\gl-$\BMGSICWY\circ\gl$-\gmres      & 1.10e-01       & 13.06       & 2          & 20          & 22          & 62          & 24          \\\hline
				\cl-$\BMGSCWY\circ\CholQR$-\gmres      & 1.02e-01       & 19.32       & 2          & 20          & 22          & 62          & 24          \\\hline
				\gl-$\BCGSPIP\circ\gl$-\gmres      & 9.96e-02       & 21.14       & 2          & 20          & 20          & 60          & 22          \\\hline
				\cl-$\BMGSICWY\circ\CholQR$-\gmres      & 9.77e-02       & 22.67       & 2          & 20          & 22          & 62          & 24          \\\hline
				\cl-$\BCGSPIP\circ\CholQR$-\gmres      & 8.55e-02       & 32.29       & 2          & 20          & 20          & 60          & 22          
			\end{tabular}
		\end{center}
		\caption{Results for \texttt{circuit\_2} example. \label{tab:circuit_2}}
	\end{table}
	
	\begin{table}[htbp]
		\begin{center}
			\begin{tabular}{c|c|c|c|c|c|c|c}
				\toprule
				\thead{Configuration}	& \thead{Time\\(s)}	& \thead{\%\\Accel.}	& \thead{Cycle\\Ct.}	& \thead{Iter.\\Ct.}	& \thead{$A$\\Ct.} & \thead{$V$\\Ct.}	& \thead{Sync\\Ct.} \\
				\midrule
				\gl-$\BCGSIROLS\circ\gl$-\gmres      & 9.54e-02       & 0.00       & 1          & 9          & 10          & 37          & 11          \\\hline
				\gl-$\BMGSICWY\circ\gl$-\gmres      & 8.55e-02       & 10.41       & 1          & 9          & 10          & 28          & 11          \\\hline
				\gl-$\BMGSCWY\circ\gl$-\gmres      & 8.11e-02       & 14.96       & 1          & 9          & 10          & 28          & 11          \\\hline
				\gl-$\BMGS\circ\gl$-\gmres      & 7.67e-02       & 19.58       & 1          & 9          & 9          & 18          & 55          \\\hline
				\cl-$\BMGSICWY\circ\CholQR$-\gmres      & 7.16e-02       & 24.91       & 1          & 8          & 9          & 25          & 10          \\\hline
				\cl-$\BMGS\circ\CholQR$-\gmres      & 6.91e-02       & 27.60       & 1          & 8          & 8          & 16          & 45          \\\hline
				\cl-$\BCGSIROLS\circ\CholQR$-\gmres      & 6.66e-02       & 30.21       & 1          & 8          & 9          & 33          & 10          \\\hline
				\gl-$\BCGSPIP\circ\gl$-\gmres      & 6.09e-02       & 36.19       & 1          & 9          & 9          & 27          & 10          \\\hline
				\cl-$\BMGSCWY\circ\CholQR$-\gmres      & 5.26e-02       & 44.92       & 1          & 8          & 9          & 25          & 10          
			\end{tabular}
		\end{center}
		\caption{Results for \texttt{rajat03} example. \label{tab:rajat03}}
	\end{table}
	
	\begin{table}[htbp!]
		\begin{center}
			\begin{tabular}{c|c|c|c|c|c|c|c}
				\toprule
				\thead{Configuration}	& \thead{Time\\(s)}	& \thead{\%\\Accel.}	& \thead{Cycle\\Ct.}	& \thead{Iter.\\Ct.}	& \thead{$A$\\Ct.} & \thead{$V$\\Ct.}	& \thead{Sync\\Ct.} \\
				\midrule
				\gl-$\BCGSIROLS\circ\gl$-\gmres      & 7.24e-01       & 0.00       & 6          & 55          & 61          & 226          & 67          \\\hline
				\gl-$\BMGSICWY\circ\gl$-\gmres      & 6.65e-01       & 8.22       & 6          & 55          & 61          & 171          & 67          \\\hline
				\gl-$\BMGSCWY\circ\gl$-\gmres      & 6.61e-01       & 8.73       & 6          & 55          & 61          & 171          & 67          \\\hline
				\gl-$\BMGS\circ\gl$-\gmres      & 6.21e-01       & 14.21       & 6          & 55          & 55          & 110          & 351          \\\hline
				\gl-$\BCGSPIP\circ\gl$-\gmres      & 5.29e-01       & 26.93       & 6          & 55          & 55          & 165          & 61          \\\hline
				\cl-$\BMGS\circ\CholQR$-\gmres      & 2.36e-01       & 67.47       & 2          & 19          & 19          & 38          & 121          \\\hline
				\cl-$\BCGSIROLS\circ\CholQR$-\gmres      & 2.02e-01       & 72.13       & 2          & 19          & 21          & 78          & 23          \\\hline
				\cl-$\BMGSCWY\circ\CholQR$-\gmres      & 1.89e-01       & 73.95       & 2          & 19          & 21          & 59          & 23          \\\hline
				\cl-$\BMGSICWY\circ\CholQR$-\gmres      & 1.86e-01       & 74.27       & 2          & 19          & 21          & 59          & 23          \\\hline
				\cl-$\BCGSPIP\circ\CholQR$-\gmres      & 1.60e-01       & 77.84       & 2          & 19          & 19          & 57          & 21          
			\end{tabular}
		\end{center}
		\caption{Results for \texttt{Kaufhold} example. \label{tab:kaufhold}}
	\end{table}
	
	\begin{table}[htbp!]
		\begin{center}
			\begin{tabular}{c|c|c|c|c|c|c|c}
				\toprule
				\thead{Configuration}	& \thead{Time\\(s)}	& \thead{\%\\Accel.}	& \thead{Cycle\\Ct.}	& \thead{Iter.\\Ct.}	& \thead{$A$\\Ct.} & \thead{$V$\\Ct.}	& \thead{Sync\\Ct.} \\
				\midrule
				\gl-$\BCGSIROLS\circ\gl$-\gmres      & 2.02e+00       & 0.00       & 15          & 144          & 159          & 591          & 174          \\\hline
				\cl-$\BMGS\circ\CholQR$-\gmres      & 1.86e+00       & 7.51       & 14          & 139          & 139          & 278          & 913          \\\hline
				\gl-$\BMGSICWY\circ\gl$-\gmres      & 1.84e+00       & 8.53       & 15          & 144          & 159          & 447          & 174          \\\hline
				\gl-$\BMGSCWY\circ\gl$-\gmres      & 1.84e+00       & 8.83       & 15          & 144          & 159          & 447          & 174          \\\hline
				\gl-$\BMGS\circ\gl$-\gmres      & 1.71e+00       & 15.17       & 15          & 144          & 144          & 288          & 939          \\\hline
				\cl-$\BCGSIROLS\circ\CholQR$-\gmres      & 1.65e+00       & 18.09       & 14          & 139          & 153          & 570          & 167          \\\hline
				\cl-$\BMGSICWY\circ\CholQR$-\gmres      & 1.54e+00       & 23.67       & 14          & 139          & 153          & 431          & 167          \\\hline
				\cl-$\BMGSCWY\circ\CholQR$-\gmres      & 1.49e+00       & 25.86       & 14          & 139          & 153          & 431          & 167          \\\hline
				\gl-$\BCGSPIP\circ\gl$-\gmres      & 1.43e+00       & 28.92       & 15          & 144          & 144          & 432          & 159          \\\hline
				\cl-$\BCGSPIP\circ\CholQR$-\gmres      & 1.26e+00       & 37.64       & 14          & 139          & 139          & 417          & 153          
			\end{tabular}
		\end{center}
		\caption{Results for \texttt{t2d\_q9} example. \label{tab:t2d_q9}}
	\end{table}
	
	\begin{table}[htbp!]
		\begin{center}
			\begin{tabular}{c|c|c|c|c|c|c|c}
				\toprule
				\thead{Configuration}	& \thead{Time\\(s)}	& \thead{\%\\Accel.}	& \thead{Cycle\\Ct.}	& \thead{Iter.\\Ct.}	& \thead{$A$\\Ct.} & \thead{$V$\\Ct.}	& \thead{Sync\\Ct.} \\
				\midrule
				\gl-$\BCGSIROLS\circ \gl$-\fom      & 5.34e+01       & 0.00       & 47          & 1162          & 1209          & 4695          & 1256          \\\hline
				\gl-$\BMGSICWY\circ \gl$-\fom      & 4.60e+01       & 13.92       & 47          & 1162          & 1209          & 3533          & 1256          \\\hline
				\gl-$\BMGSCWY\circ \gl$-\fom      & 4.57e+01       & 14.41       & 47          & 1162          & 1209          & 3533          & 1256          \\\hline
				\cl-$\BMGSSVL\circ \MGSSVL$-\fom      & 4.39e+01       & 17.77       & 45          & 1121          & 1121          & 3363          & 34890          \\\hline
				\cl-$\BMGSLTS\circ \MGSLTS$-\fom      & 4.29e+01       & 19.63       & 45          & 1121          & 1121          & 3363          & 34890          \\\hline
				\gl-$\BCGSPIP\circ \gl$-\fom      & 4.07e+01       & 23.70       & 47          & 1162          & 1162          & 3486          & 1209          \\\hline
				\cl-$\BCGSIROLS\circ \CholQR$-\fom      & 4.06e+01       & 23.91       & 45          & 1121          & 1166          & 4529          & 1211          \\\hline
				\cl-$\BMGS\circ \CholQR$-\fom      & 3.97e+01       & 25.60       & 45          & 1121          & 1121          & 2242          & 15697          \\\hline
				\gl-$\BMGSSVL\circ \gl$-\fom      & 3.94e+01       & 26.21       & 47          & 1162          & 1162          & 3486          & 3533          \\\hline
				\gl-$\BMGS\circ \gl$-\fom      & 3.89e+01       & 27.15       & 47          & 1162          & 1162          & 2324          & 16237          \\\hline
				\gl-$\BMGSLTS\circ \gl$-\fom      & 3.88e+01       & 27.24       & 47          & 1162          & 1162          & 3486          & 3533          \\\hline
				\cl-$\BMGSCWY\circ \CholQR$-\fom      & 3.72e+01       & 30.25       & 45          & 1121          & 1166          & 3408          & 1211          \\\hline
				\cl-$\BMGSICWY\circ \CholQR$-\fom      & 3.66e+01       & 31.44       & 45          & 1121          & 1166          & 3408          & 1211          \\\hline
				\cl-$\BCGSPIP\circ \CholQR$-\fom      & 3.41e+01       & 36.12       & 181          & 2162          & 2162          & 6486          & 2343          
			\end{tabular}
		\end{center}
		\caption{Results for \lapl example. \label{tab:lapl_2d}}
	\end{table}
\end{appendix}
  
\end{document}